\documentclass[a4paper]{article}
\usepackage{geometry}
\usepackage[utf8]{inputenc}
\usepackage{graphicx}
\usepackage{amsmath}
\usepackage{amssymb}
\usepackage{float}
\usepackage{caption}
\usepackage{hyperref}
\usepackage{algorithm}
\usepackage{xcolor}
\usepackage{mathtools}

\usepackage[noend]{algpseudocode}

\usepackage[backend=bibtex, style=alphabetic]{biblatex}
\def \trans{^{\scriptscriptstyle{\intercal}}}

\def \ep{\hbox{ }\hfill$\Box$}

\mathtoolsset{showonlyrefs}

\DeclareMathOperator{\R}{\mathbb{R}}

\def \E{\mathbb{E}}
\def \F{\mathbb{F}}
\def \M{\mathbb{M}}
\def \S{\mathbb{S}}

\def \P{\mathbb{P}}

\def\1{{\bf 1}}
\def \I{{\bf I}}
\def \N{\mathbb{N}}
\newcommand{\un}{1\hspace{-1mm}{\rm I}}   %
\def\Fc{{\cal F}}
\def\Tc{{\cal T}}
\def\Uc{{\cal U}}

\def\Wc{{\cal W}}
\def\Xc{{\cal X}}
\def\Zc{{\cal Z}}
\def\Nc{{\cal N}}
\def\argmin_#1{\underset{#1}{\mathrm{argmin\, }}}
\newtheorem{remark}{Remark}
\newcommand{\di}{\mathrm{d}}

\def \Sum{\displaystyle\sum}

\numberwithin{equation}{section} 
\numberwithin{remark}{section}

\def \mrx{\mathrm{x}}

\begin{document}
\title{Neural networks-based backward scheme for fully nonlinear PDEs
\thanks{This work is supported by  FiME, Laboratoire de Finance des March\'es de l'Energie, and the ''Finance and Sustainable Development'' EDF - CACIB Chair.}}
\author{
Huy\^en \sc{Pham}\footnote{LPSM, Universit\'e de Paris,  CREST-ENSAE \& FiME \sf \href{mailto:pham at lpsm.paris}{pham at lpsm.paris}} 
\and
Xavier \sc{Warin}
\footnote{EDF R\&D \& FiME \sf \href{mailto:xavier.warin at edf.fr}{xavier.warin at edf.fr}} 
\and
Maximilien \sc{Germain}
\footnote{EDF R\&D, LPSM, Université de Paris  \sf \href{mailto:maximilien.germain at edf.fr}{mgermain at lpsm.paris}}
}

\date{This version: December 7, 2020}

\maketitle

\begin{abstract}
We propose a numerical method for solving high dimensional fully nonlinear partial differential equations (PDEs).  Our algorithm estimates simultaneously 
by backward time induction the solution and its gradient by multi-layer neural networks, 
while the Hessian is approximated by automatic differentiation of the gradient at previous step. 
This methodology extends to the fully non\-linear case 
the approach recently proposed in \cite{HPW19} for semi-linear PDEs. Numerical tests  illustrate the performance and accuracy of our method on several examples 
in high dimension with non-linearity on the Hessian term  including a linear quadratic control problem with control on the diffusion coefficient, 
Monge-Amp\`ere equation and Hamilton-Jacobi-Bellman equation in portfolio optimization. 
\end{abstract}

\vspace{5mm}

\noindent {\bf Key words:} Neural networks,  fully nonlinear PDEs in high dimension, backward scheme.

\vspace{5mm}

\noindent {\bf MSC Classification:}  60H35, 65C20, 65M12.

\section{Introduction}

This paper is devoted to the resolution in high dimension  of fully nonlinear parabolic partial differential equations (PDEs) of the form 
\begin{equation} \label{eq:PDEInit}
\left\{
\begin{aligned}
 \partial_t u + f(.,.,u,D_xu,D_x^2u) & = 0 , \;\;\;\;\;\; \mbox{ on } [ 0,T)\times\R^d, \\
 u(T,.) &=g, \;\;\;\;\;  \mbox{ on } \R^d,
 \end{aligned}
 \right.
 \end{equation}
 with a  non-linearity  in the solution, its gradient $D_x u$  and its hessian $D_x^2 u$ via the function $f(t,x,y,z, \gamma)$ defined on $[0,T]\times\R^d\times\R\times\R^d\times\S^d$ (where $\S^d$ is the set of symmetric  $d\times d$ matrices), and a terminal condition $g$. 
  
 The numerical resolution of this class of PDEs is far more difficult than the one of classical semi-linear PDEs  where the nonlinear function $f$ does not depend on $\gamma$. 
 In fact, rather few methods are available to solve fully nonlinear equations even in moderate dimension.
 \begin{itemize}
 \item First based on the work of
 \cite{cheridito2007second}, an effective scheme  developed in \cite{fahim2011probabilistic} using some regression techniques  has been shown to be convergent under some ellipticity  conditions later removed by \cite{tan2013splitting}. Due to the use of basis functions, this scheme does not permit to solve PDE in dimension greater than 5.
 \item A scheme based on nesting Monte Carlo has been recently proposed in \cite{warin2018monte}. It seems to be effective in very high dimension for maturities $T$ not too long and  linearities not too important. 
 \item A numerical algorithm to solve fully nonlinear equations has been proposed by \cite{BEJ19} based on the second order backward stochastic differential equations (2BSDE) representation of \cite{cheridito2007second} and  global deep neural networks minimizing a terminal objective function, but no test  on real fully nonlinear case is given.  
 This extends the idea introduced in the pioneering papers \cite{weinan2017deep,han2017solving}, which were the first serious works 
 for using machine learning methods to solve high dimensional PDEs. 
 \item The Deep Galerkin method proposed in \cite{sirignano2018dgm} based on some machine learning techniques and using some automatic  differentiation of the solution seems to be effective on some cases. It has been tested in \cite{al2018solving} for example on the Merton problem.
 \end{itemize}
 
 In this article,  we introduce a numerical method  based on machine learning techniques  and backward in time  iterations, which extends the proposed schemes in \cite{szpruch19} for linear problems, and in the recent work   \cite{HPW19} for  semi-linear PDEs.  
 The approach in these works consists in estimating  simultaneously  the solution and its gradient by multi-layer neural networks by minimizing a sequence of loss functions defined in backward induction.  
 A basic idea to extend this method to the fully nonlinear case would rely on  the representation proposed in \cite{cheridito2007second}:  
 at each time step  $t_n$ of an Euler scheme, the Hessian $D_x^2u$ at $t_n$ is approximated by a neural network minimizing 
 some local $L_2$ criterion associated to a BSDE involving $D_x u$ at date $t_{n+1}$ and $D_x^2 u$. Then, the pair $(u,D_x u)$ at date $t_n$ is approximated/learned with a second minimization 
 similarly as in the method described by \cite{HPW19}. The first minimization can be implemented with different variations but numerical results show that the global scheme  does 
 not scale well with the dimension. Instability on the  $D_x^2u$ calculation rapidly propagates during the backward resolution.  Besides, the methodology appears to be costly when using two optimizations at each time step.   
 An alternative approach that we develop here, is to combine the ideas of \cite{HPW19}  and the splitting method in \cite{beck2019deep} in order  to derive a new deep learning scheme that requires only  one local optimization during the backward resolution for learning the pair $(u,D_x u)$ and approximating $D_x^2 u$ by automatic differentiation of the gradient computed  at the previous step. 

The outline of the paper is organized as follows.  In Section \ref{secnum}, we briefly recall the mathematical description of the classical feedforward approximation, and then derive 
the proposed neural networks-based backward scheme. We test our method in Section \ref{sectest} on various examples. 
First we illustrate our results with  a PDE  involving a non-linearity of type $u D_x^2u$.  Then, we consider a stochastic linear quadratic problem with controlled volatility where an analytic solution is available, and we test the performance and accuracy of our algorithm up to dimension $20$. Next, we apply our algorithm to a Monge-Amp\`ere equation, and finally, we provide numerical tests for the solution to fully nonlinear Hamilton-Jacobi-Bellman equation, with non-linearities of the form $|D_x u|^2/D_x^2 u$,  arising in portfolio selection problem with stochastic volatilities. 

\section{The proposed deep backward scheme} \label{secnum}

Our aim is to numerically approximate the function $u$ $:$ $[0,T]\times\R^d$ $\mapsto$ $\R$, assumed to be the unique smooth solution to the fully nonlinear PDE \eqref{eq:PDEInit} under suitable conditions.  This will be achieved by means of neural networks approximations for $u$ and its gradient $D_x u$, relying on a backward scheme and training simulated data of some forward diffusion process.  Approximations of PDE in high dimension by neural networks have now become quite popular, and are supported theoretically by recent results in \cite{hutzenthaler2018overcoming} and \cite{DLM19} showing their efficiency to overcome the curse of dimensionality.

 \subsection{Feedforward neural network to approximate functions}

We denote by $d_0$ the dimension of  the input variables, and $d_1$ the dimension of the output variable.  
A (deep) neural  network is characterized by a number of layers  $L+1$ $\in$ $\N\setminus\{1,2\}$  with $m_\ell$, $\ell$ $=$ $0,\ldots,L$, the number of neurons (units or nodes) on each layer: the first layer is the input layer with $m_0$ $=$ $d$, the last layer is the output layer with $m_L$ $=$ $d_1$, and the $L-1$ layers between are called hidden layers, where we choose for simplicity the same  dimension $m_\ell$ $=$ $m$, $\ell$ $=$ $1,\ldots,L-1$.
 
A  feedforward neural network is a  function from  $\R^{d_0}$ to $\R^{d_1}$ defined as the composition
\begin{align} \label{defNN}
x \in \R^d  & \longmapsto  \; A_L \circ  \varrho \circ A_{L - 1} \circ \ldots \circ \varrho \circ A_1(x) \; \in \; \R. 
\end{align}
Here $A_\ell$, $\ell$ $=$ $1,\ldots,L$ are affine transformations: $A_1$ maps from $\R^{d_0}$ to $\R^m$, $A_2,\ldots,A_{L-1}$ map from $\R^m$ to $\R^m$, and $A_L$ maps from $\R^m$ to 
$\R^{d_1}$, represented by 
\begin{align}
A_\ell (x) &= \; \Wc_\ell x + \beta_\ell,
\end{align}
for a matrix $\Wc_\ell$ called weight, and a vector $\beta_\ell$ called  bias term,  $\varrho$ $:$ $\R$ $\rightarrow$ $\R$ is a nonlinear function, called activation function, and applied 
component-wise on the outputs of $A_\ell$, i.e., $\varrho(x_1,\ldots,x_m)$ $=$ $(\varrho(x_1),\ldots,\varrho(x_m))$. Standard examples of activation functions are the sigmoid, the ReLu, the Elu, $\tanh$. 

All these matrices $\Wc_\ell$ and vectors $\beta_\ell$, $\ell$ $=$ 
$1,\ldots,L$,  are the parameters of the neural network, and can be identified with  an element $\theta$ $\in$ $\R^{N_m}$, where $N_m$ $=$ 
$\sum_{\ell=0}^{L-1} m_\ell (1+m_{\ell+1})$ $=$ $d_0(1+m)+m(1+m)(L-2)+m(1+d_1)$ is the number of parameters.  
We denote by  $\Nc_{d_0,d_1,L,m}$ the set of all functions generated by  \eqref{defNN} for $\theta \in \R^{N_m}$.

\subsection{Forward-backward representation}

Let us introduce a forward diffusion process 
\begin{align}
 X_t  &= X_0 +   \int_0^t \mu(s,X_s) ds+  \int_0^t \sigma(s,X_s)  dW_s, \;\;\; 0 \leq t\leq T, 
 \label{eq:SDE}
\end{align}
where $\mu$ is a function defined on $[0,T] \times \R^d$ with values in $\R^d$, $\sigma$ is a  function defined on $[0,T] \times \R^d$ with values in  $\M^d$ the set of  $d \times d$ matrices, and 
$W$ a $d$-dimensional Brownian motion  on some probability space $(\Omega,\Fc,\P)$ equipped with a filtration  $\F$ $=$ $(\Fc_t)_{0\leq t\leq T}$ satisfying the usual conditions. 
The process $X$ will be used for the simulation of training data in our deep learning algorithm, and we shall discuss later the choice of the drift and diffusion coefficients $\mu$ and $\sigma$, see Remark \ref{remchoiceX}. 

Let us next denote by $(Y,Z,\Gamma)$ the triple of $\F$-adapted processes valued in $\R\times\R^d\times\S^d$,  defined by
\begin{align} \label{relYU}
Y_t  \;= \;  u(t,X_t), \;\;\;  Z_t \; = \; D_x u(t,X_t), \;\;\;   \Gamma_t \: = \; D_x^2 u(t,X_t), \;\;\; 0 \leq t \leq T. 
\end{align}
By It\^o's formula applied to $u(t,X_t)$, and since $u$ is solution to  \eqref{eq:PDEInit},  we see that $(Y,Z,\Gamma)$ satisfies the backward equation: 
\begin{align}
Y_t & = \; g(X_T) -  \int_t^T \big[ \mu(s,X_s).Z_s + \frac{1}{2}{\rm tr}(\sigma\sigma\trans(s,X_s) \Gamma_s) - f(s,X_s,Y_s,Z_s,\Gamma_s) \big] ds \nonumber  \\
& \;\;\;\;\;\;\;\;\;\;   - \int_t^T \sigma\trans(s,X_s) Z_s. dW_s, \;\;\;\;\;  0 \leq t \leq T.  \label{BSDEY}
\end{align}

\begin{remark}
{\rm 
This BSDE does not uniquely characterize a triple $(Y,Z,\Gamma$) contrarily to the semilinear case (without a non-linearity with respect to $\Gamma$) in which proper assumptions on the equation coefficients provide existence and uniqueness for a solution couple $(Y,Z)$. In the present case at least two options can be used to estimate the $\Gamma$ component:
\begin{itemize}
    \item Rely on the 2BSDE representation from \cite{cheridito2007second} which extends the probabilistic representation of \cite{pardoux1990adapted} for semilinear equations to the fully nonlinear case. It is the approach used by \cite{BEJ19} with a global large minimization problem, as in \cite{han2017solving}.
    \item Compute the second order derivative by automatic differentiation. This is the point of view we adopt in this paper together with a local approach solving several small optimization problems. In this way,  we provide an extension of \cite{HPW19} to cover a broader range of nonlinear PDEs.
\end{itemize}
}
\ep
\end{remark}

\subsection{Algorithm}

We now provide a numerical approximation of the forward backward system \eqref{eq:SDE}-\eqref{BSDEY}, and consequently of the solution $u$ (as well as its gradient $D_x u$) to the PDE 
\eqref{eq:PDEInit}. 

We start from a time grid $\pi$ $=$ $\{t_i, i=0,\ldots,N\}$ of $[0,T]$, with $t_0$ $=$ $0$ $<$ $t_1$ $<$  $\ldots$ $<$ $t_N$ $=$ $T$, and time steps $\Delta t_i$ $:=$ $t_{i+1}-t_i$, $i$ $=$ $0,\ldots,N-1$.  
The time discretization of the forward process $X$ on $\pi$  is then equal (typically when $\mu$ and $\sigma$ are constants) or approximated by an Euler scheme: 
\begin{align*}
X_{t_{i+1}} & =  \; X_{t_i} + \mu(t_i,X_{t_i}) \Delta t_i + \sigma(t_i,X_{t_i}) \Delta W_{t_i}, \;\;\;\;\; i =0,\ldots,N-1, 
\end{align*}
where we set $\Delta W_{t_i}$ $:=$ $W_{t_{i+1}}-W_{t_i}$ (by misuse of notation, we keep the same notation $X$ for the continuous time diffusion process and its Euler scheme).  The backward SDE \eqref{BSDEY} is approximated by the time discretized scheme
\begin{align*}
Y_{t_i} & \simeq  \; Y_{t_{i+1}}  -  \big[ \mu(t_i,X_{t_i}). Z_{t_i} + \frac{1}{2} {\rm tr}\big(\sigma\sigma\trans(t_i,X_{t_i} \Gamma_{t_i}\big) 
- f(t_i,X_{t_i},Y_{t_i},Z_{t_i},\Gamma_{t_i}) \big] \Delta t_i - \sigma\trans(t_i,X_{t_i})Z_{t_i}.\Delta W_{t_i}, 
\end{align*}
that is written in forward form as
\begin{align} \label{Yforward} 
Y_{t_{i+1}} & \simeq  \; F(t_i,X_{t_i},Y_{t_i},Z_{t_i},\Gamma_{t_i},\Delta t_i,\Delta W_{t_i}),  \;\;\; i = 0, \ldots,N-1,  
\end{align}
with 
\begin{align} 
F(t,x,y,z,\gamma,h,\Delta)  & := \;  y  - \tilde f(t,x,y,z,\gamma) h  \; + \;   z\trans\sigma (t,x)  \Delta, \label{expressF} \\
\tilde f(t,x,y,z,\gamma) & := \; f(t,x,y,z,\gamma) -  \mu(t,x).z -   \frac{1}{2} {\rm tr}\big( \sigma \sigma\trans(t,x) \gamma\big).  \nonumber
\end{align}

The idea of the proposed scheme is the following. 
Similarly as in \cite{HPW19}, we approximate at each time $t_i$, $u(t_i,.)$ and its gradient $D_x u(t_i,.)$, by neural networks $x$ $\in$ $\R^d$ $\mapsto$ 
$(\Uc_i(x;\theta),\Zc_i(x;\theta))$ with parameter $\theta$ that are learned optimally by backward induction: suppose that $\hat \Uc_{i+1}$ $:=$ $\Uc_{i+1}(.;\theta_{i+1}^*)$, 
$\hat\Zc_{i+1}$ $:=$ $\Zc_{i+1}(.;\theta_{i+1}^*)$ is an approximation of $u(t_{i+1},.)$ and $D_x u(t_{i+1},.)$ at time $t_{i+1}$, then $\theta_i^*$ is computed from the minimization of the quadratic loss function: 
\begin{align}
\hat L_i(\theta) &= \; \E \Big| \hat \Uc_{i+1}  -    F(t_i,X_{t_i},\Uc_i(X_{t_i};\theta),\Zc_i(X_{t_i};\theta), D\hat\Zc_{i+1}(\Tc(X_{t_{i+1}})),\Delta t_i,\Delta W_{t_i}) \Big|^2
\end{align}
where  $\Tc$ is a truncation operator such that $\Tc(X)$ is bounded for example by a quantile of the diffusion process and $D\hat\Zc_{i+1}$  
stands for the automatic  differentiation of $\hat\Zc_{i+1}$.  The idea behind the truncation is the following. During  one step resolution,  the estimation of the gradient is less accurate at the edge  of the explored domain where samples are rarely generated. Differentiating the gradient gives a very oscillating Hessian at the edge of the domain. At the following time step resolution, these oscillations propagate to the gradient and the solution even if the domain where the oscillations occur is  rarely attained. In order to avoid these oscillations,  a truncation is achieved, 
permits to  avoid that the  oscillations of the neural network fit  in zone where the simulations  propagate scarcely  to areas of importance.  This truncation may be necessary to get convergence on some rather difficult cases.  
Of course this truncation is only valid if the real Hessian does not varies too much.

The intuition for the relevance of this scheme to the approximation of the PDE \eqref{eq:PDEInit} is the following. From \eqref{relYU} and \eqref{Yforward}, the solution $u$ to \eqref{eq:PDEInit} should 
approximately satisfy
\begin{align*}
u(t_{i+1},X_{t_{i+1}}) & \simeq \;  F(t_i,X_{t_i},u({t_i},X_{t_i}),D_x u({t_i},X_{t_i}),D_x^2 u({t_i},X_{t_i}),\Delta t_i,\Delta W_{t_i}).
\end{align*}
Suppose that at time $t_{i+1}$, $\hat\Uc_{i+1}$ is an estimation of $u(t_{i+1},.)$. Recalling 
the expression of $F$ in \eqref{expressF}, the quadratic loss function at time $t_i$ is then approximately equal to
\begin{align*}
\hat L_i(\theta) & \simeq \;  \E \Big| u(t_{i},X_{t_{i}}) -  \Uc_i(X_{t_i};\theta)  + \big( D_x u(t_i,X_{t_i}) - \Zc_i(X_{t_i};\theta) \big)\trans \sigma(t_i,X_{t_i})\Delta W_{t_i} \\
& \;\;\; - \;  \Delta t_i \big[ \tilde f(t_i,X_{t_i},u({t_i},X_{t_i}),D_x u({t_i},X_{t_i}),D_x^2 u({t_i},X_{t_i})) 
- \tilde f(t_i,X_{t_i},\Uc_i(X_{t_i};\theta),\Zc_i(X_{t_i};\theta), D\hat\Zc_{i+1}(\Tc(X_{t_{i+1}}))) \big]  \Big|^2. 
\end{align*}
By assuming that $\tilde f$ has small non-linearities in its arguments $(y,z,\gamma)$, say Lipschitz, possibly with a suitable choice of $\mu$, $\sigma$, the loss function is thus approximately equal to 
\begin{align*}
\hat L_i(\theta) & \simeq \;  (1 + O(\Delta t_i)) \E\big| u(t_{i},X_{t_{i}}) -  \Uc_i(X_{t_i};\theta) \big|^2 + O(\Delta t_i) \E\big| D_x u(t_{i},X_{t_{i}}) -  \Zc_i(X_{t_i};\theta) \big|^2 + O(|\Delta t_i|^2).  
\end{align*}
Therefore, by minimizing over $\theta$ this quadratic loss function, via stochastic gradient descent (SGD) based on simulations of $(X_{t_i},X_{t_{i+1}},\Delta W_{t_i})$ (called training data in the machine learning language), one expects the neural networks $\Uc_i$ and $\Zc_i$ to learn/approximate better and better the functions $u(t_i,.)$ and $D_x u(t_i,)$ in view of the universal approximation theorem for neural networks.  
The rigorous convergence of this algorithm is postponed to a future work.

To sum up, the global algorithm is given in Algo \ref{theAlgorithm} in the case where $g$ is Lipschitz and the derivative can be analytically calculated almost everywhere.  
If the derivative of $g$ is not available, it can be calculated by automatic differentiation of the neural network approximation of $g$. 

\begin{algorithm}[H]
\caption{\label{theAlgorithm} Algorithm for fully non-linear equations.}
\begin{algorithmic}[1]
 \State  Use a single deep neural network $(\Uc_N(.;\theta),\Zc_N(.;\theta))$ $\in$ $\Nc_{d,1+d,L,m}$ and  minimize (by SGD)
 \begin{equation} \label{eq:initScheme}
     \left\{ 
     \begin{aligned} 
     \hat L_N(\theta)  & := \;  \E  \Big| \Uc_{N}(X_{t_{N}};\theta) - g(X_{t_{N}}) \Big|^2  + \frac{\Delta t_{N-1}}{d}  \E  \Big| \Zc_{N}(X_{t_{N}};\theta) - Dg(X_{t_{N}}) \Big|^2 \\
     \theta_N^* & \in {\rm arg}\min_{\theta\in\R^{N_m}} \hat L_N(\theta). 
     \end{aligned}
     \right.
     \end{equation}
 \State  $\widehat\Uc_N$ $=$ $\Uc_N(.;\theta_N^*)$, and set $\widehat\Zc_N$ $=$ $\Zc_N(.;\theta_N^*)$
 \For{ $i$ $ =$ $N-1,\ldots,0$}
 \State Use a single deep neural network $(\Uc_i(.;\theta),\Zc_i(.;\theta))$ $\in$ $\Nc_{d,1+d,L,m}$
 for the approximation of $(u(t_i,.),D_x u(t_i,.))$, and compute (by SGD) the minimizer of the expected quadratic loss function
     \begin{equation} \label{eq:scheme}
     \left\{ 
     \begin{aligned} 
     \hat L_i(\theta)  & := \;  \E \Big| \widehat\Uc_{i+1}(X_{t_{i+1}}) \; - \;   
     F(t_i,X_{t_i},\Uc_i(X_{t_i};\theta),\Zc_i(X_{t_i};\theta),D \hat \Zc_{i+1}(\Tc(X_{t_{i+1}})),\Delta t_{i},\Delta W_{t_{i}}) \Big|^2  \\
     \theta_i^* & \in {\rm arg}\min_{\theta\in\R^{N_m}} \hat L_i(\theta). 
     \end{aligned}
     \right.
     \end{equation}
 \State  Update: $\widehat\Uc_i$ $=$ $\Uc_i(.;\theta_i^*)$, and set $\widehat\Zc_i$ $=$ $\Zc_i(.;\theta_i^*)$.
 \EndFor
 \end{algorithmic}
\end{algorithm}

\begin{remark}\label{rem: implicit explicit}
{\rm
Several alternatives can be implemented for the computation of the second order derivative. A natural candidate would consist in choosing to approximate the solution $u$ at time $t_i$ by a neural network $\Uc_i$ and estimate $\Gamma_i$ as the iterated automatic differentiation $D_x^2 \Uc _i$. However, it is shown in  \cite{HPW19} that choosing only a single neural network for $u$ and using its automatic derivative to estimate the $Z$ component degrades the error in comparison to the choice of two neural networks $\Uc,\Zc$. A similar behavior has been observed during our tests for this second order case and the most efficient choice was to compute the derivative of the $\Zc$ network. 
This derivative can also be estimated at the current time step $t_i$ instead of $t_{i+1}$. However this method  leads to an additional cost for the neural networks training by complicating the computation of the automatic gradients performed by Tensorflow during the backpropagation. It also leads numerically to worse results on the control estimation, as empirically observed in Table \ref{fig: table results Merton Implicit} and described in the related paragraph "Comparison with an implicit version of the scheme".  For this reason, we decided  to apply a splitting method and evaluate the Hessian at time $t_{i+1}$.  For this reason, we decided  to apply a splitting method and evaluate the Hessian at time $t_{i+1}$. 
\ep}
\end{remark}

\begin{remark} \label{remchoiceX} 
{\rm The diffusion process $X$ is used for the training simulations in the stochastic gradient descent method for finding the minimizer of the quadratic loss function in \eqref{eq:scheme}, 
where the expectation is replaced by empirical average for numerical implementation.  The choice of the drift and diffusion parameters are explained in Section \ref{secchoice}.  
}
\ep
\end{remark}


\section{Numerical results} \label{sectest}

We first  construct an example  with different non-linearities in the Hessian term and  the solution.  We graphically show that the solution is very well calculated in dimension $d=1$ and then move to higher dimensions. We then use an example derived from a stochastic optimization problem with an analytic solution and show that we are able to accurately calculate the solution.  Next, we consider the numerical resolution of the Monge-Amp\`ere equation, and finally, give some tests for a fully nonlinear Hamilton-Jacobi-Bellman equation arising from portfolio optimization with stochastic volatilities.

\subsection{Choice of the algorithm hyperparameters}  \label{secchoice} 


We describe in this paragraph how we choose the various hyperparameters of the algorithm and explain the learning strategy. 

\vspace{1mm}

\noindent $\bullet$  {\sc Parameters of the training simulations}:
the choice of the drift coefficient  is typically  related to the underlying probabilistic problem associated to the PDE (for example a stochastic control problem), and should drive the training process  to regions of interest, e.g.., that are visited with large probability by the optimal state process in stochastic control. In practice, we can take a drift function $\mu(.)$  equal to the drift associated to some  a priori  control. This choice of control could be an optimal control for a related problem for which we know the solution, or could be the control obtained by the first iteration of the algorithm. 
The choice of the diffusion coefficient $\sigma$ is also important:  large $\sigma$  induces a better exploration of the state space,  
but as  we will see in most of examples below,  it gives a scheme slowly converging to the solution with respect to the time discretization and it generates a higher variance on the results.
Moreover, for the applications in stochastic control, we might explore some region that are visited with very small probabilities by the optimal state process, hence representing few interest. On the other hand, small $\sigma$ means a weak exploration, and we might lack  information and precision on some region of the state space:  
the solution calculated at each time step is far more sensitive to very local errors induced by the neural network approximation and tends to generate a bias. Therefore a trade off has to be found between rather high variance with slow convergence  in time  and  fast convergence in time with a potential bias. We also refer to \cite{nusric20} for a discussion on the role of the diffusion coefficient. 

In practice and for the numerical examples in the next section, we test the scheme for different $\sigma$ and by varying the number of time steps, and if it converges to the same solution, one can consider  that we have obtained the correct solution. 
We also  show the impact of the choice of the diffusion coefficient $\sigma$. 

\vspace{1mm}

\noindent $\bullet$  {\sc Parameters of truncation}:
Given the training simulations $X$, we choose a truncation operator $\Tc_p$ indexed by a parameter $p$  close to $1$, so that  $\Tc_p(X_t^{})$ corresponds to a truncation of $X_t^{}$ at a given quantile $\phi_p$.  
In the numerical tests, we shall vary $p$ between $0.95$ and $0.999$. 

\vspace{1mm}

\noindent $\bullet$ {\sc Parameters of the optimization algorithm over neural networks}:
In the whole numerical part, we use a classical Feedforward network using layers with $m$ neurons each and a $\tanh$ activation function, the output layer uses an identity activation function.  At each time step the resolution of equation \eqref{eq:scheme} is achieved using a mini-batch with $1000$ training trajectories. The training and learning rate adaptation procedure is the following:
\begin{itemize}
    \item Every 40  inner gradient descent iterations, the loss is checked on $10000$ validation trajectories. 
    \item This optimization sequence is repeated with $200$ outer iterations for the first optimization step at date $t_N=T$ and only $100$ outer iterations  at the dates $t_i$ with $i <N$. 
    \item An average  of the loss calculated on $10$ successive outer iterations is performed. If the decrease of the average loss every $10$ outer iterations is  less than $5\%$ then the learning rate is divided by $2$. 
\end{itemize}

The optimization is performed using the Adam  gradient descent algorithm, see \cite{kingma2014adam}.  Notice that the adaptation of the learning rate is not common with the Adam method but in our case it appears to be crucial to have a steady diminution of the loss of the objective function. The procedure  is also described in \cite{chan2019machine} and the chosen parameters are similar to this article. At the initial optimization step at time $t_N=T$, the learning rate is taken equal to $1E-2$ and at the following optimization steps,  we start with a learning rate equal to $1E-3$.


During time resolution, it is far more effective to initialize the solution of equations \eqref{eq:scheme}  with the solution $(\Uc,\Zc)$ at the next time step. Indeed the previously computed values at time step $t_{i+1}$ are 
good approximations of the processes at time step $t_{i}$ if the PDE solution and its gradient are continuous.
All experiments are achieved using Tensorflow  \cite{2015tensorflow}.
In the sequel, the PDE solutions on curves are calculated as the average of 10 runs. We provide the standard deviation associated to these results. We also  show the influence of the number of neurons on the accuracy of the results.


\subsection{A non-linearity in $u D_x^2 u$}

 We consider a generator in the form
 \begin{align*}
 f(t, x,y,z,\gamma) &= \;  
 y {\rm tr}( \gamma ) + \frac{y}{2} + 2 y^2 -2 y^4 e^{-(T-t)}, 
 \end{align*}
 and $g(x)$ $=$ $\tanh{\big(\frac{\sum_{i=1}^d x_i}{\sqrt{d}}\big)}$, so  that an analytical solution is available:
 \begin{align*}
 u(t, x) &= \;  \tanh{\Big(\frac{\sum_{i=1}^d x_i}{\sqrt{d}}\Big)} e^{-\frac{T-t}{2}}.
 \end{align*}
 
We fix the horizon $T$ $=$ $1$, and choose to evaluate the solution at $t=0$ and  $x =0.5 \frac{\un_d}{\sqrt{d}}$ (here $\un_d$ denotes the vector in $\R^d$ with all components equal to $1$), 
 for which $u(t,x)$ $=$ $0.761902$ while its derivative is equal to $1.2966$. 
 This initial value $x$ is chosen such that independently of the dimension the solution is varying around this point and not in a region where the $\tanh$ function is close  to $-1$ or $1$. 
 
The coefficients of the forward process used to solve the equation are  (here $\I_d$ is the identity $d\times d$-matrix) 
 \begin{flalign*}
\sigma = & \frac{\hat \sigma}{\sqrt{d}} \I_d,  \quad  \mu  = 0,  
\end{flalign*}
and here the truncation operator  is chosen equal to
\begin{align*}
 \Tc_p(X_t^{0,x}) = \min\big\{\max[x  -  \sigma \sqrt{t} \phi_p ,  X_t^{0,x} ],  x+  \sigma \sqrt{t} \phi_p \big\},
\end{align*}
where $\phi_p$ $=$  $\mathcal{N}^{-1}(p)$,  with  $\mathcal{N}$ is the CDF of a unit centered Gaussian random variable. 

\begin{figure}[h!]
\begin{minipage}[b]{0.32\linewidth}
  \centering
 \includegraphics[width=\textwidth]{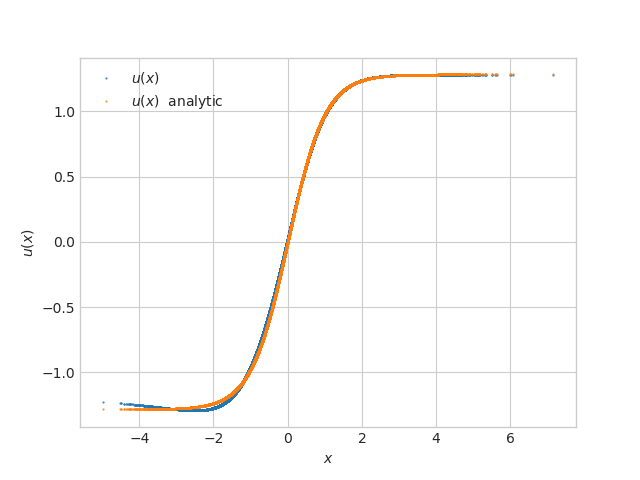}
 \caption*{$Y$ at date $t=0.5$.}
 \end{minipage}
\begin{minipage}[b]{0.32\linewidth}
  \centering
 \includegraphics[width=\textwidth]{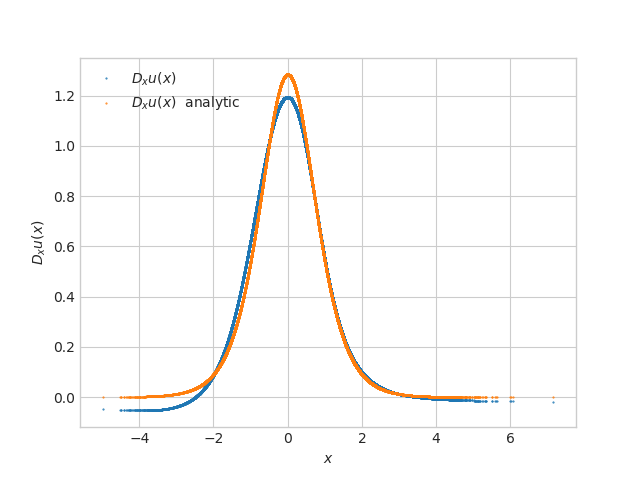}
 \caption*{$Z$ at date $t=0.5$ }
 \end{minipage}
 \begin{minipage}[b]{0.32\linewidth}
  \centering
 \includegraphics[width=\textwidth]{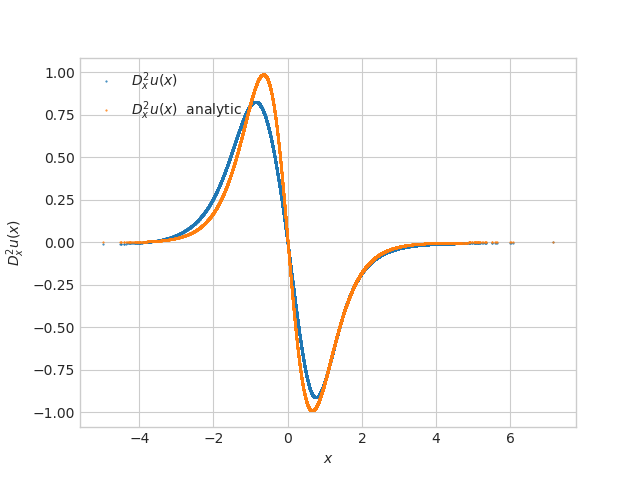}
 \caption*{$\Gamma$ at date $t=0.5$ }
 \end{minipage}
 \begin{minipage}[b]{0.32\linewidth}
  \centering
 \includegraphics[width=\textwidth]{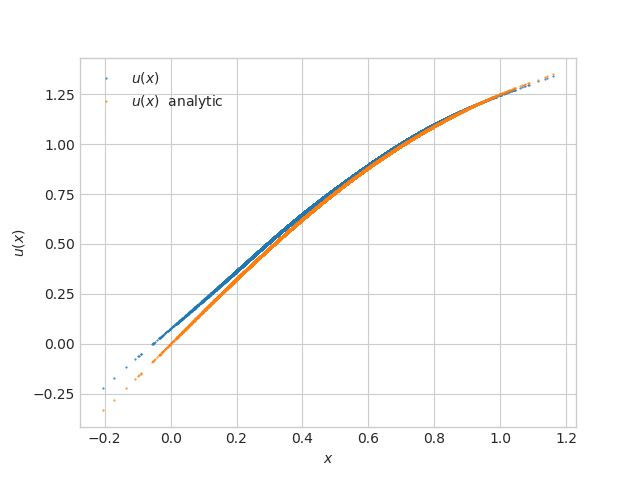}
 \caption*{$Y$ at date $t=0.006125$.}
 \end{minipage}
\begin{minipage}[b]{0.32\linewidth}
  \centering
 \includegraphics[width=\textwidth]{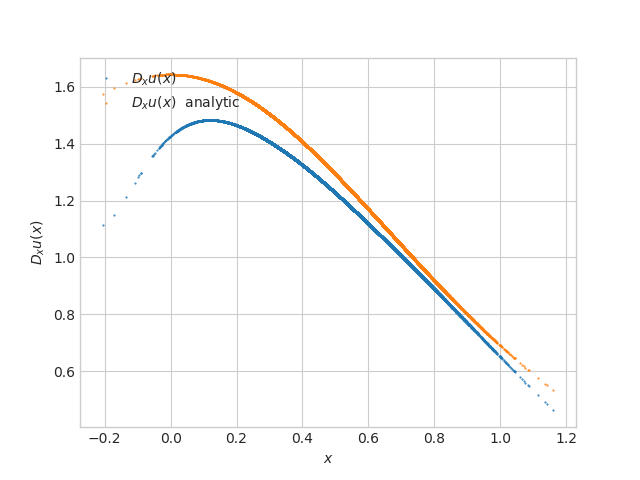}
 \caption*{$Z$ at date $t=0.006125$ }
 \end{minipage}
 \begin{minipage}[b]{0.32\linewidth}
  \centering
 \includegraphics[width=\textwidth]{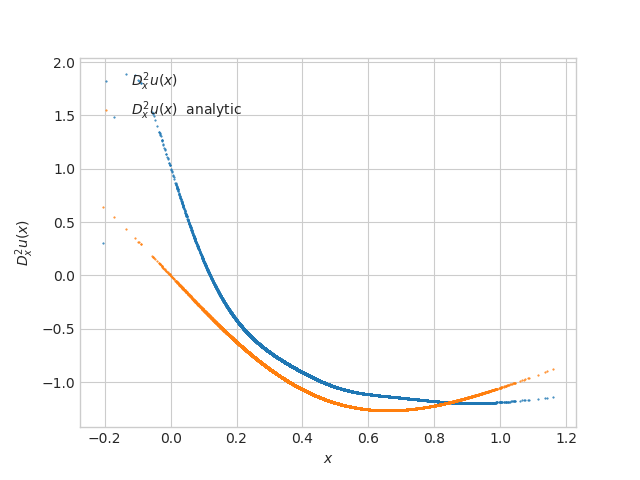}
 \caption*{$\Gamma$ at date $t=0.006125$ }
 \end{minipage}
 \caption{\label{fig:linLap1D}  A single valuation run for test case one  1D  using $160$ time steps, $\hat \sigma=2.$, $p =0.999$, 20 neurons, 2 layers.}
 \end{figure}

In the numerical results,  we take $p=0.999$ and $m$ $=$  $20$ neurons. 
We first begin in dimension $d=1$, and show in  Figure \ref{fig:linLap1D}  how  $u$, $D_x u$ and $D_x^2u$ are well approximated by the resolution method. 

\begin{figure}[H]
\begin{minipage}[b]{0.49\linewidth}
  \centering
 \includegraphics[width=\textwidth]{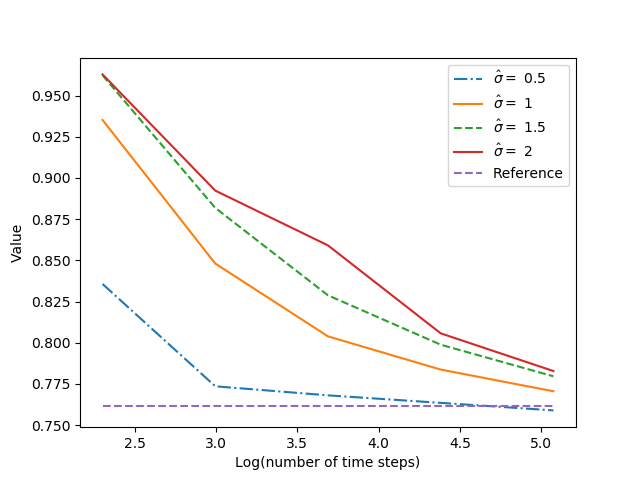}
 \caption*{Convergence of $u$  depending on  $\hat \sigma$}
 \end{minipage}
\begin{minipage}[b]{0.49\linewidth}
  \centering
 \includegraphics[width=\textwidth]{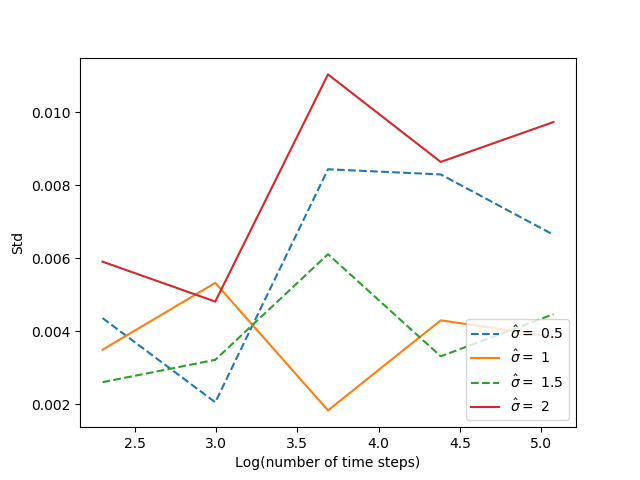}
 \caption*{Standard   deviation of $u$}
 \end{minipage}
 \begin{minipage}[b]{0.49\linewidth}
  \centering
 \includegraphics[width=\textwidth]{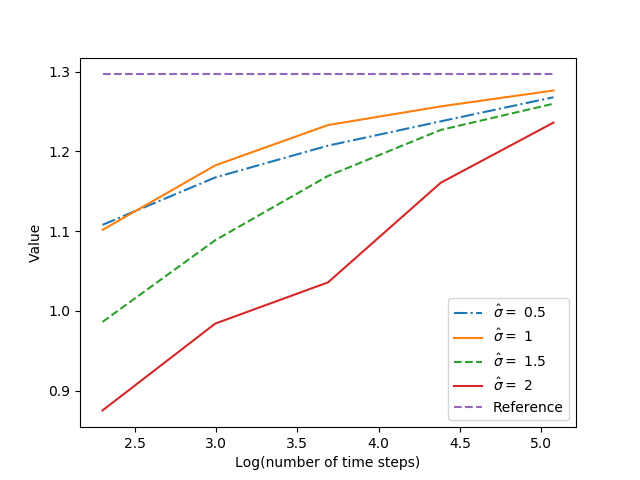}
 \caption*{Convergence of $D_xu$  depending on  $\hat \sigma$}
 \end{minipage}
\begin{minipage}[b]{0.49\linewidth}
  \centering
 \includegraphics[width=\textwidth]{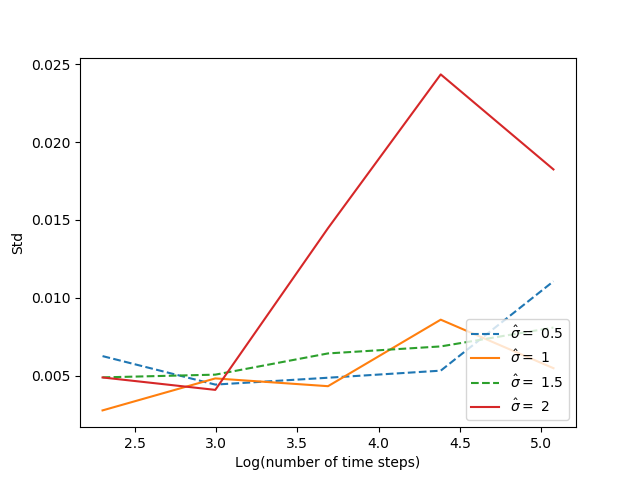}
 \caption*{Standard    deviation of $D_x u$}
 \end{minipage}
 \caption{\label{fig:linearLap1D} Convergence in 1D of the case one, number of neurons par layer equal to $20$, 2 layers, $p=0.999$.}
 \end{figure}

 On Figure \ref{fig:linearLap1D}, we check the convergence, for different values of $\hat \sigma$ of both the solution $u$  and its derivative at point $x$ and date $0$. Standard deviation of the function value is very low and the standard deviation of the derivative  still being low.

\begin{figure}[H]
\begin{minipage}[b]{0.49\linewidth}
  \centering
 \includegraphics[width=\textwidth]{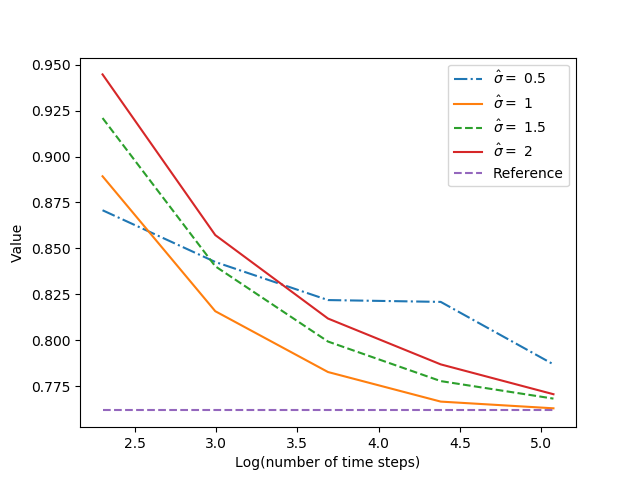}
 \caption*{Convergence of $u$  depending on  $\hat \sigma$}
 \end{minipage}
 \begin{minipage}[b]{0.49\linewidth}
  \centering
 \includegraphics[width=\textwidth]{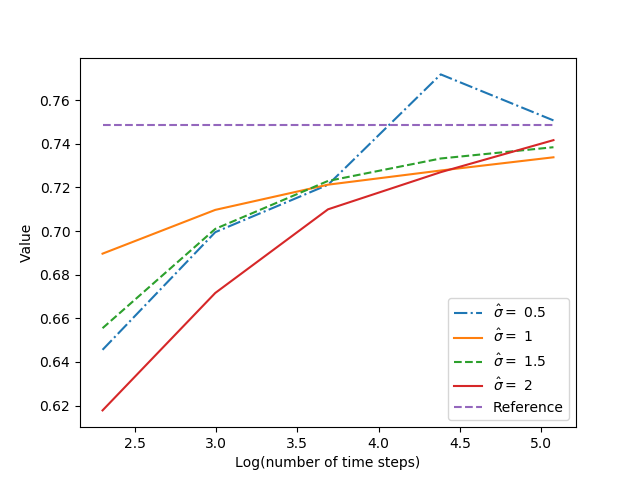}
 \caption*{\small{Convergence of  $D_x u$ (first component) depending on  $\hat \sigma$}}
 \end{minipage}
 \caption{\label{fig:linearLap3D} Convergence in 3D of the case one, number of neurons par layer equal to $20$, 2 layers, $p=0.999$.}
 \end{figure}

\begin{figure}[H]
\begin{minipage}[b]{0.49\linewidth}
  \centering
 \includegraphics[width=\textwidth]{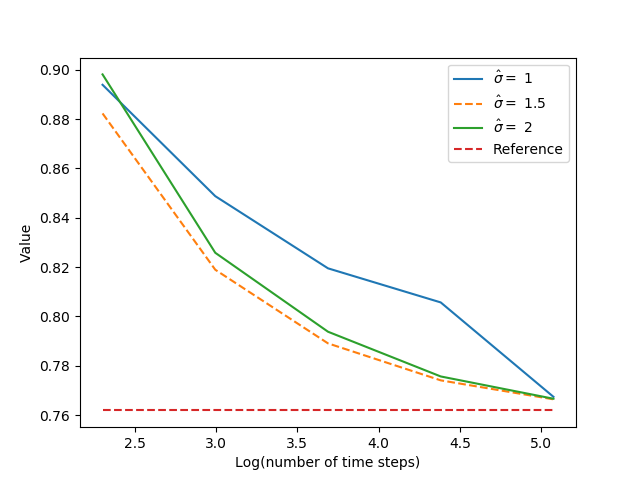}
 \caption*{Convergence of $u$  depending on  $\hat \sigma$}
 \end{minipage}
 \begin{minipage}[b]{0.49\linewidth}
  \centering
 \includegraphics[width=\textwidth]{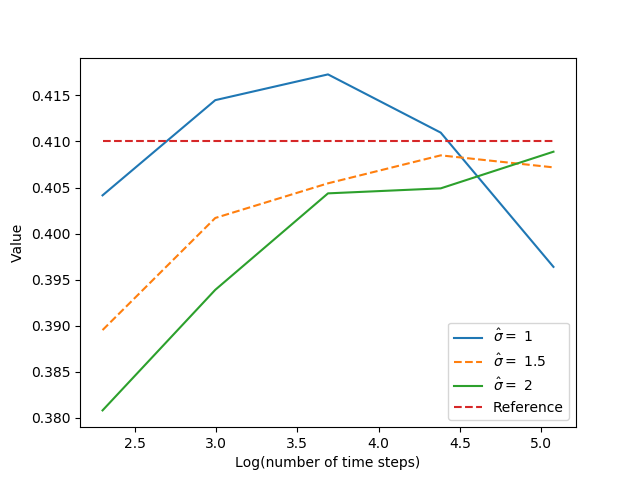}
 \caption*{\small{Convergence of  $D_xu$  depending on  $\hat \sigma$ (first component)}}
 \end{minipage}
 \caption{\label{fig:linearLap10D} Convergence in 10D of the case one, number of neurons par layer equal to $20$, 2 layers, $p=0.999$.}
 \end{figure}

 As the dimension increases, we have to increase the value of $\hat \sigma$  of the forward process.
In dimension 3, the value $\hat \sigma =0.5$ gives high  standard deviation in the result obtained as shown on  Figure \ref{fig:linearLap3D}, while in dimension 10, see Figure \ref{fig:linearLap10D}, we see that the value $\hat \sigma=1$  is too low to give good results. We also clearly notice that in 10D, a smaller time step should be used but in our test cases we decided to consider a maximum number of time steps equal to 160.

 On this simple test case, the dimension is not a problem and very good results are obtained in dimension 20 or above with only 20 neurons and 2 layers.

 \clearpage

 \subsection{A linear quadratic stochastic test case.}

 In this example,  we consider a controlled process $\Xc$ $=$ $\Xc^\alpha$ with dynamics in $\R^d$ according to 
 \begin{align*}
  d\Xc_t  &= \;  (A \Xc_t+ B \alpha_t) dt + D\alpha_t dW_t, \;\;\; 0 \leq t \leq T, \;\;   \Xc_0  \; = \;   x, 
 \end{align*}
 where $W$ is a real Brownian motion, the control process  $\alpha$ is  valued in $\R$, and the constant  coefficients $A \in \M^d$, $B \in \R^{d}$, $D \in \R^d$.
 The quadratic cost functional to be minimized is 
 \begin{flalign*}
  J(\alpha) =& \E\Big[ \int_0^T \big( \Xc_t\trans Q \Xc_t + \alpha_t^2 N \big) dt + \Xc_T\trans P \Xc_T  \Big],
 \end{flalign*}
 where $P$, $Q$ are non negative $d \times d$ symmetric matrices and $N \in \R$ is strictly positive. 
 
 The Bellman equation associated to this stochastic control problem is:
 \begin{align*}
  \frac{\partial u}{ \partial t } + \inf_{a \in \R} \big[  (A x + B a).D_x u + \frac{a^2}{2} {\rm tr}( DD\trans D_x^2 u) + x\trans Q x  + N a^2 \big] & =  \; 0, \;\;\; (t,x) \in [0,T)\times\R^d, \\ 
  u(T,x) &= \;  x\trans P x, \;\;\; x  \in \R^d,
 \end{align*}
which can be rewritten as a fully nonlinear equation in the form  \eqref{eq:PDEInit} with
\begin{align*}
f(t,x,y,z,\gamma) & = \; x\trans Q x +  Ax.z  - \frac{1}{2} \frac{| B\trans z|^2}{{\rm tr}(D D\trans \gamma)+2N}.   
\end{align*}
 An explicit solution to this  PDE is given by 
 \begin{align*}
  u(t,x) &= \; x\trans K(t) x,
 \end{align*}
 where $K(t)$ is non negative $d \times d$ symmetric matrix function  solution to the Riccati equation:
 \begin{flalign*}
  \dot{K} + A^\top K + K A + Q - \frac{K B B^\top K}{N + D^\top K D}= 0,  \quad  K(T)= P.
 \end{flalign*}

 We take $T=1$. The coefficients of the  forward process used to solve the equation are  
 \begin{flalign*}
\sigma = & \frac{\hat \sigma}{\sqrt{d}} \I_d,  \quad  \mu(t,x) = A x.
\end{flalign*}
In our numerical example we take the following parameters for the optimization problem: 
\begin{flalign*}
 A = \I_d, \;  B =D = \un_d,  & \quad Q =P = \frac{1}{d} \I_d, \quad  N=d
\end{flalign*} 
and we want to estimate the solution at $x=\un_d$.

In this example, the truncation operator (indexed by $p$ between $0$ and $1$ and close to $1$)  is as follows:
\begin{align*}
 \Tc_p(X_t^{x}) &= \; \min\big\{  \max\big[x e^{ \hat A t} -  \sigma \sqrt{\frac{e^{ 2\hat A t} - \hat 1}{2 \hat A}} \phi_p,  X_t^{x} \big],x e^{ \hat A t} +  \sigma \sqrt{\frac{e^{ 2\hat A t} - \hat 1}{2 \hat A}} \phi_p\big\},
\end{align*}
where $\phi_p$ $=$  $\mathcal{N}^{-1}(p)$,  $\hat A$ is a vector so that $\hat A_i= A_{ii}$, $i=1, ..., d$,  $\hat 1$ is a unit vector, and the square root is taken componentwise.

\begin{figure}[H]
\begin{minipage}[b]{0.32\linewidth}
  \centering
 \includegraphics[width=\textwidth]{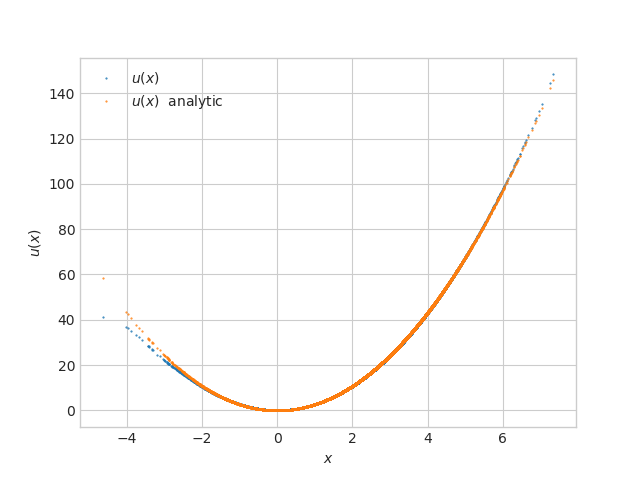}
 \caption*{$Y$ at date $t=0.5$.}
 \end{minipage}
\begin{minipage}[b]{0.32\linewidth}
  \centering
 \includegraphics[width=\textwidth]{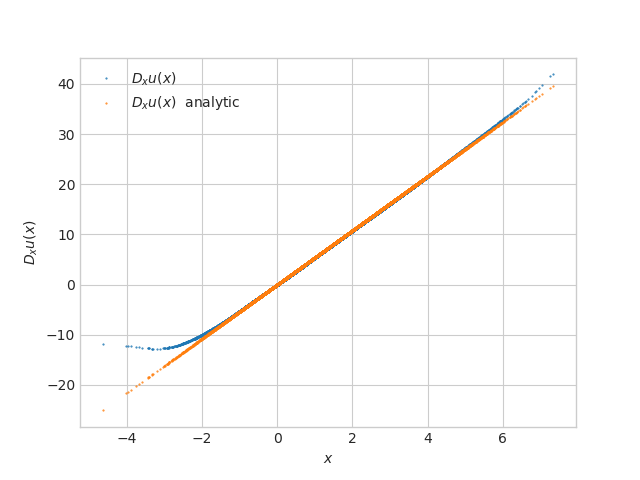}
 \caption*{$Z$ at date $t=0.5$ }
 \end{minipage}
 \begin{minipage}[b]{0.32\linewidth}
  \centering
 \includegraphics[width=\textwidth]{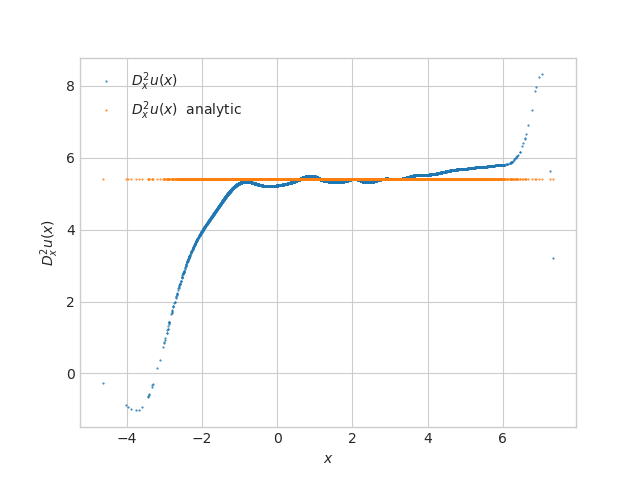}
 \caption*{$\Gamma$ at date $t=0.5$ }
 \end{minipage}
 \begin{minipage}[b]{0.32\linewidth}
  \centering
 \includegraphics[width=\textwidth]{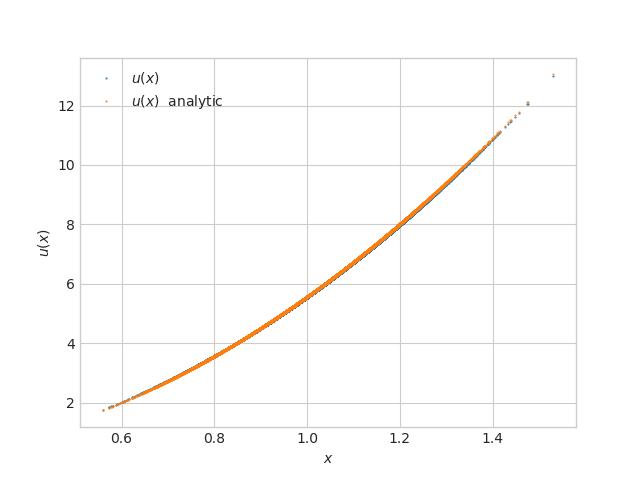}
 \caption*{$Y$ at date $t=0.006125$.}
 \end{minipage}
\begin{minipage}[b]{0.32\linewidth}
  \centering
 \includegraphics[width=\textwidth]{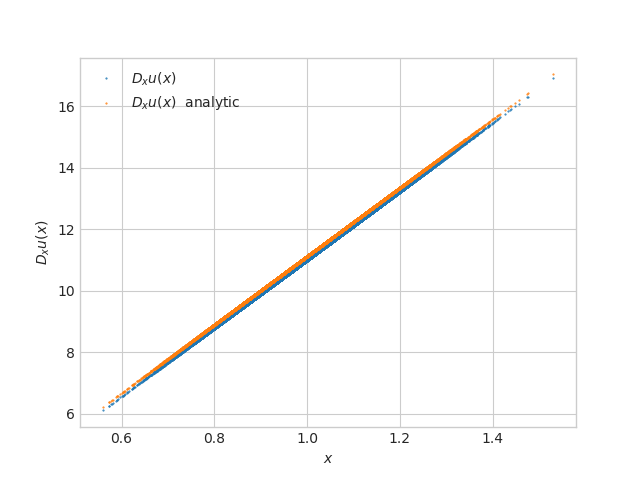}
 \caption*{$Z$ at date $t=0.006125$ }
 \end{minipage}
 \begin{minipage}[b]{0.32\linewidth}
  \centering
 \includegraphics[width=\textwidth]{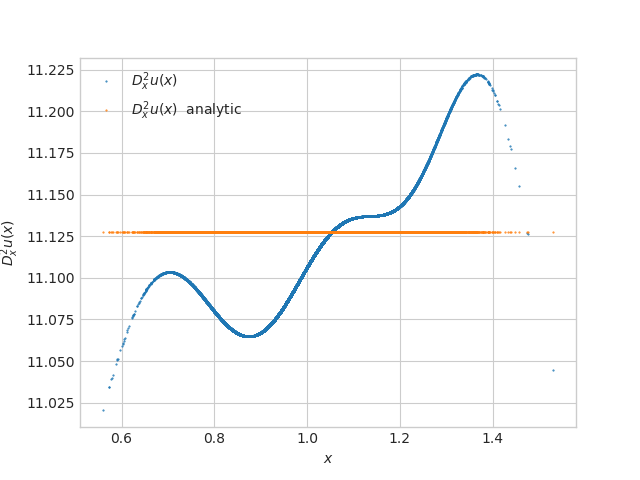}
 \caption*{$\Gamma$ at date $t=0.006125$ }
 \end{minipage}
 \caption{\label{fig:linearquad1D}  Test case linear quadratic 1D  using $160$ time steps, $\hat \sigma=1.5$, $p =0.999$, $100$ neurons.}
 \end{figure}

On Figure \ref{fig:linearquad1D} we give the solution of the PDE with $d=1$ using $\hat \sigma=1.5$ obtained for two dates: at $t=0.5$ and at $t$ close to zero. We observe that we have a very good estimation of the function value and a correct one of the $\Gamma$ value  at date $t=0.5$. The precision remains good  for $\Gamma$ close to $t=0$  and 
very good for $u$ and $D_x u$.
 
\begin{figure}[H]
\begin{minipage}[b]{0.49\linewidth}
  \centering
 \includegraphics[width=\textwidth]{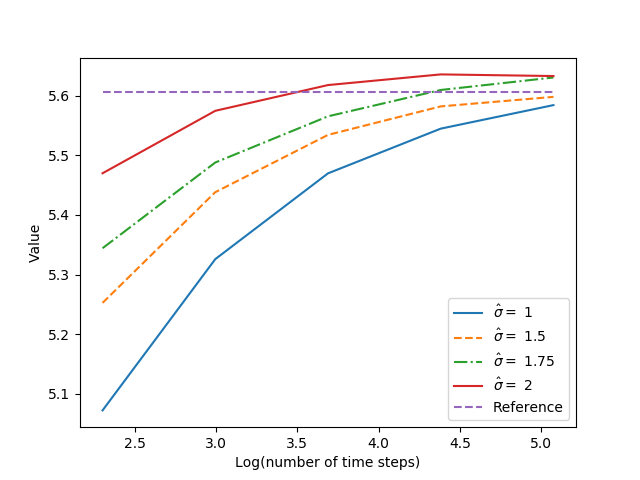}
 \caption*{Convergence   of $u$ depending on  $\hat \sigma$.}
 \end{minipage}
\begin{minipage}[b]{0.49\linewidth}
  \centering
 \includegraphics[width=\textwidth]{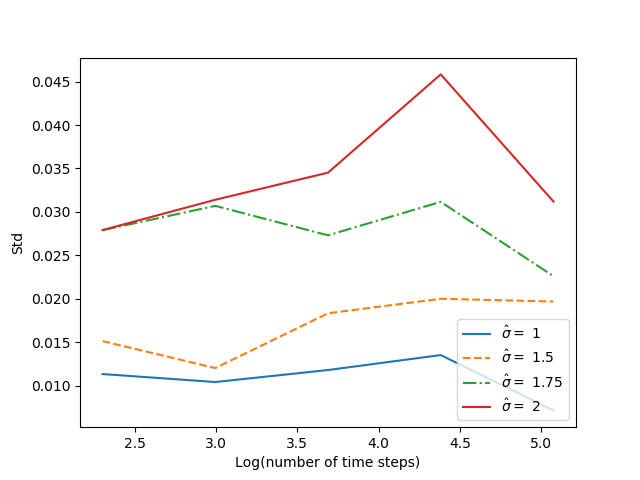}
 \caption*{Standard   deviation of $u$ }
 \end{minipage}
 \caption{\label{fig:linearquad1DB} Convergence in 1D of the linear quadratic case, number of neurons par layer equal to $50$, 2 layers, $p=0.999$.}
 \end{figure}

 On Figure \ref{fig:linearquad1DB}, we give the results obtained in dimension $d=1$  by varying  $\hat \sigma$. 
 For a value of $ \hat \sigma=2$,  the standard deviation of the result becomes far higher than with $\hat \sigma = 0.5$ or $1.$

\begin{figure}[H]
\begin{minipage}[b]{0.49\linewidth}
  \centering
 \includegraphics[width=\textwidth]{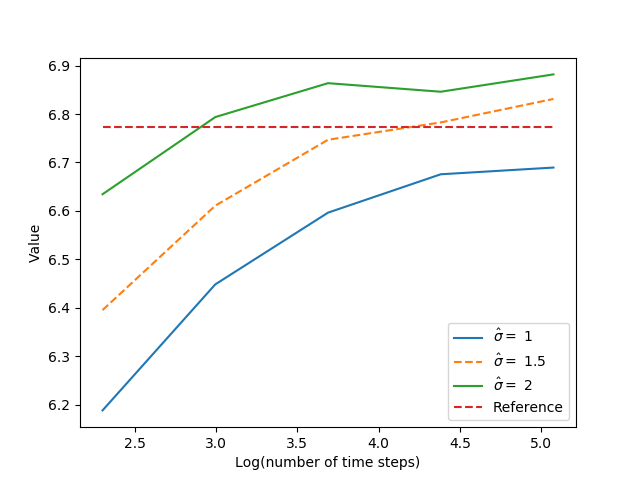}
 \caption*{ 10 neurons}
 \end{minipage}
\begin{minipage}[b]{0.49\linewidth}
  \centering
 \includegraphics[width=\textwidth]{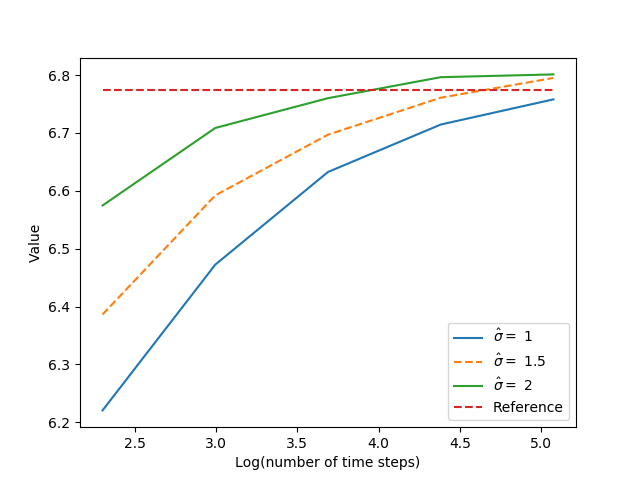}
 \caption*{ 20 neurons}
 \end{minipage}
 \begin{minipage}[b]{0.49\linewidth}
  \centering
 \includegraphics[width=\textwidth]{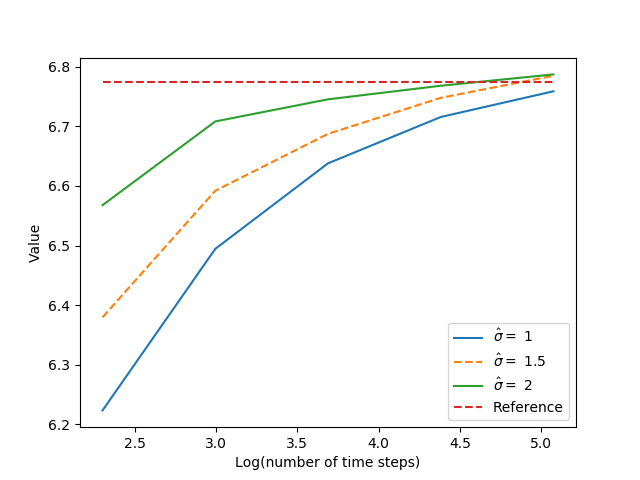}
 \caption*{ 30 neurons}
 \end{minipage}
 \begin{minipage}[b]{0.49\linewidth}
  \centering
 \includegraphics[width=\textwidth]{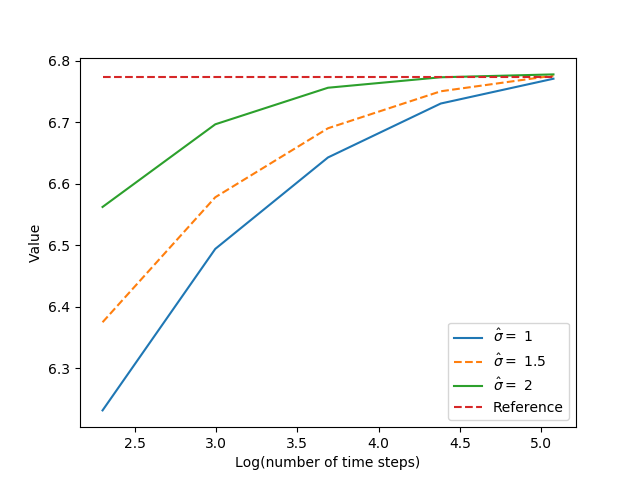}
 \caption*{ 50 neurons}
 \end{minipage}
 \caption{\label{fig:linearquad3D}  Convergence in 3D of the linear quadratic case, 2 layers, testing the influence of the number of neurons, truncation $p=0.95$.}
 \end{figure}

 On Figure \ref{fig:linearquad3D}, for $d=3$, we take a quite low truncation factor $p=0.95$ and observe that  the number of  neurons to take has to be rather high. 
 We have also checked that taking a number of hidden layers equal to 3 does not  improve the results.

\begin{figure}[H]
\begin{minipage}[b]{0.49\linewidth}
  \centering
 \includegraphics[width=\textwidth]{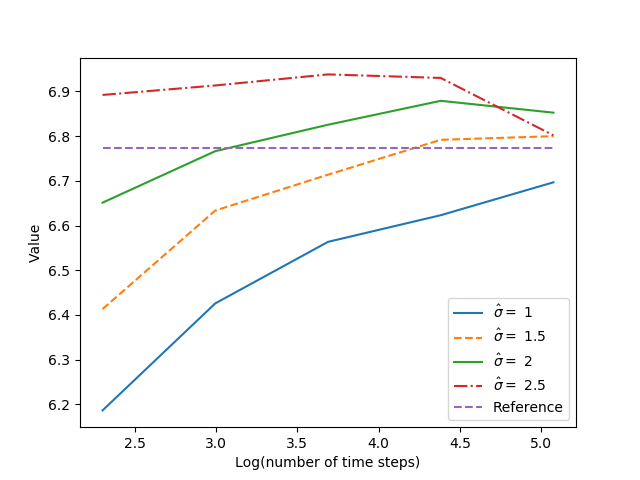}
 \caption*{ 10 neurons}
 \end{minipage}
\begin{minipage}[b]{0.49\linewidth}
  \centering
 \includegraphics[width=\textwidth]{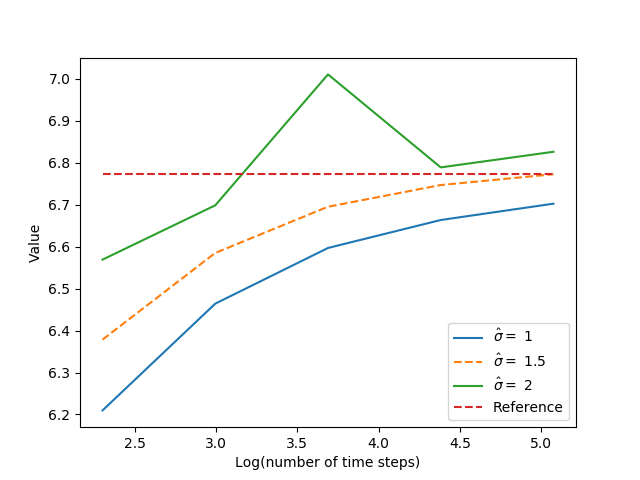}
 \caption*{ 20 neurons}
 \end{minipage}
 \begin{minipage}[b]{0.49\linewidth}
  \centering
 \includegraphics[width=\textwidth]{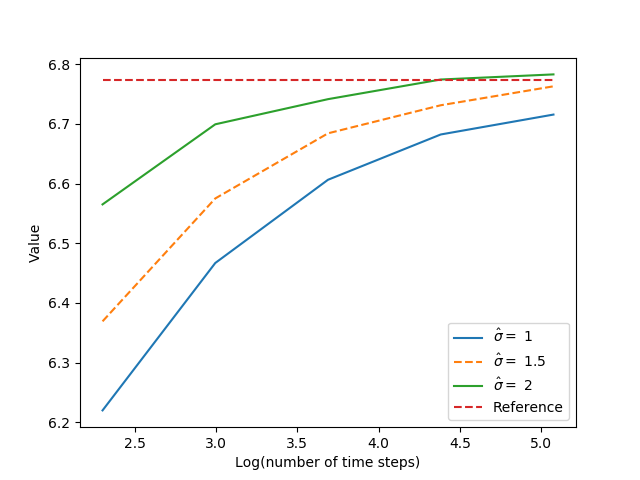}
 \caption*{ 50 neurons}
 \end{minipage}
 \begin{minipage}[b]{0.49\linewidth}
  \centering
 \includegraphics[width=\textwidth]{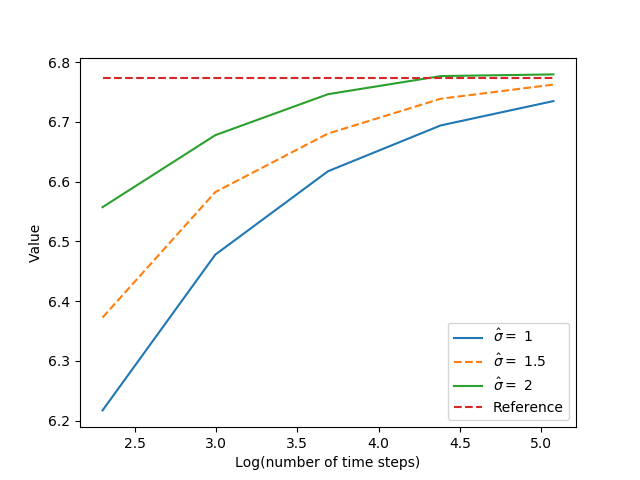}
 \caption*{ 100 neurons}
 \end{minipage}
 \caption{\label{fig:linearquad3D99}  Convergence in 3D of the linear quadratic case, 2 layers, testing the influence of the number of neurons, truncation $p=0.99$.}
 \end{figure}

 On Figure \ref{fig:linearquad3D99}, for $d=3$, we give the same graphs for a higher truncation factor. 
 As we take a higher truncation factor,  the results are improved by taking a  higher number  of neurons ($100$ in the figure below).
  

\begin{figure}[H]
 \begin{minipage}[b]{0.48\linewidth}
  \centering
 \includegraphics[width=\textwidth]{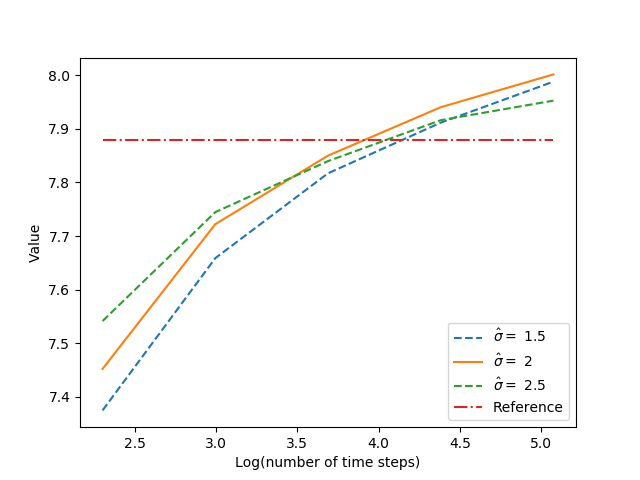}
 \caption*{Convergence with  50 neurons}
 \end{minipage}
  \begin{minipage}[b]{0.48\linewidth}
  \centering
 \includegraphics[width=\textwidth]{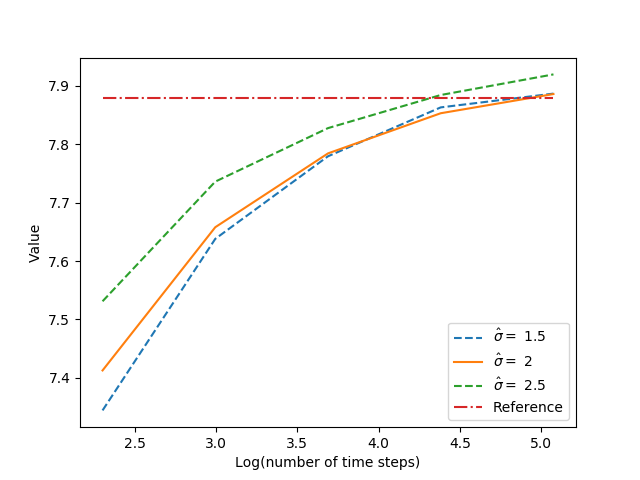}
 \caption*{ Convergence with 100 neurons}
 \end{minipage}
 \caption{\label{fig:linearquad7D}  Convergence in 7D of the linear quadratic case, 2 layers, $p=0.999$.}
 \end{figure}

 On Figure \ref{fig:linearquad7D}, we observe in dimension 7 the influence of the number of neurons on the result for a high truncation factor $p=0.999$.
 We clearly have a bias for a number of neurons equal to $50$. This bias disappears when the number of neurons increases to $100$.
 
\begin{figure}[H]
  \centering
 \includegraphics[width=0.48\textwidth]{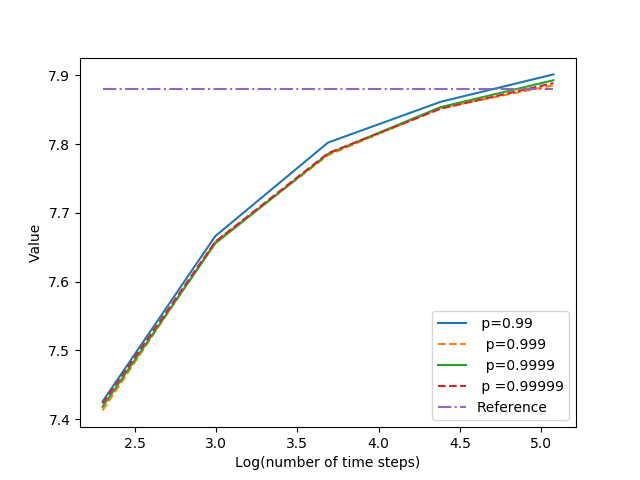}
 \caption{\label{fig:linearquad7DB} Function value convergence in 7D of the linear quadratic case with   2 layers, 100 neurons, testing $p$, using $\hat \sigma=2$}.
 \end{figure}

 On Figure \ref{fig:linearquad7DB}, for $d=7$, we check that influence of the truncation factor  appears to be slow  for higher dimensions.
 
\begin{figure}[H]
 \begin{minipage}[b]{0.49\linewidth}
  \centering
 \includegraphics[width=\textwidth]{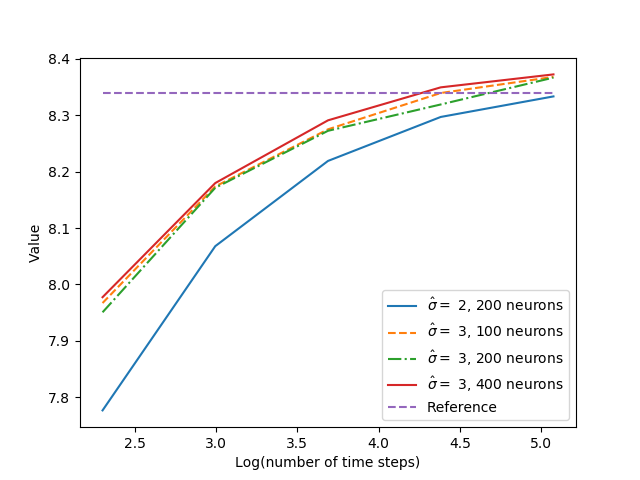}
 \caption*{ 10D}
 \end{minipage}
 \begin{minipage}[b]{0.49\linewidth}
  \centering
 \includegraphics[width=\textwidth]{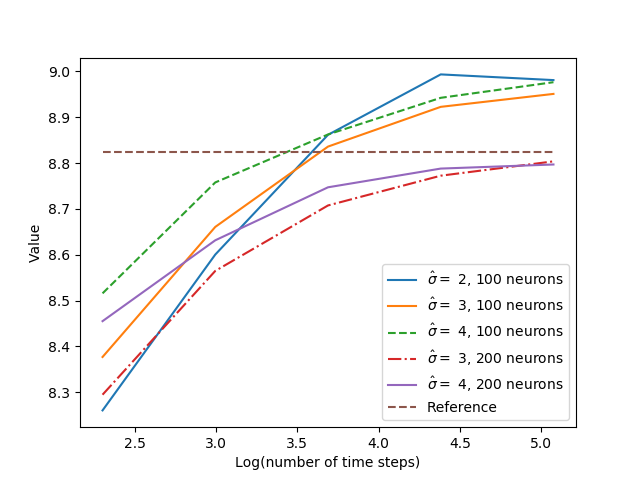}
 \caption*{ 15D}
 \end{minipage}
 \caption{\label{fig:linearquad1015D}  Function value convergence in 10D and 15D of the linear quadratic case with  2 layers, $p=0.999$.}
 \end{figure}

\begin{figure}[H]
  \centering
 \includegraphics[width=0.48\textwidth]{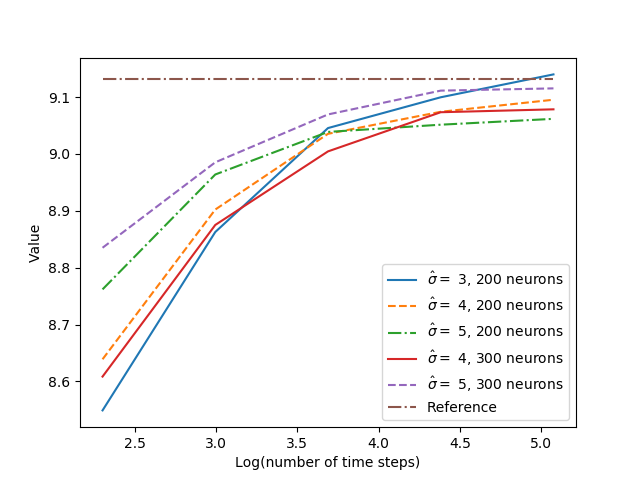}
 \caption{\label{fig:linearquad20D}  Function value convergence in 20D of the linear quadratic case with  2 layers, $p=0.999$.}
 \end{figure}

Finally, we give results in dimension 10, 15 and 20 for $p=0.999$ on Figures \ref{fig:linearquad1015D}, \ref{fig:linearquad20D}. We observe that the number a neurons with 2 hidden layers has to increase with the dimension but also that the increase is rather slow in contrast with the case of one hidden layer as theoretically shown in  \cite{pinkus1999approximation}.
For $\hat \sigma =5$ we had to take 300 neurons to get very accurate results.

\subsection{Monge-Amp\`ere equation}     


Let us consider the  parabolic Monge-Amp\`ere equation 
\begin{align} \label{Monge} 
\begin{cases}
\partial_t u + {\rm det}(D_x^2 u) \; = \; h(x), \quad  (t,x)  \in [0,T]\times \R^d, \\
u(T,x) \; = \; g(x), 
\end{cases}
\end{align}
 where ${\rm det}(D_x^2 u)$ is the determinant of the Hessian matrix $D_x^2 u$. It is in the form \eqref{eq:PDEInit} with 
\begin{align}
f(t,x,\gamma) &= \; {\rm det}(\gamma) - h(x).
\end{align}

 \vspace{1mm}

We  test our algorithm by choosing a  $C^2$ function $g$, then compute 
$G$ $=$  ${\rm det}( D_x^2 g)$, and set $h$ $:=$ $G-1$.   Then, by construction, the function
\begin{align}
u(t,x) &= \; g(x) + T-t, 
\end{align}
is solution to the Monge-Amp\`ere  equation \eqref{Monge}.  We choose $g(x) = \cos(\sum_{i=1}^d x_i/\sqrt{d})$,  
and we shall train with  the forward process  $X$ $=$ $x_0 + W$, where $W$ is a $d$-dimensional Brownian motion. On this example, 
we use neural networks with 3 hidden layers, $d+10$ neurons per layer,  
and we do not need to apply any truncation to the forward process $X$. Actually, we observe that adding a truncation worsens the results.
For choosing the truncation level, we first test the method with no truncation before decreasing the quantile parameter $p$. In the Monge-Ampère case the best results are obtained without any truncation. It may be caused by the oscillation of the Hessian. 

\vspace{2mm}

The following table gives the results in dimension $d$ $=$ $5$, $15$, and for $T$ $=$ $1$.

\vspace{2mm}

\begin{table}[H]
    \centering
   \begin{tabular}{|c|c|c|c|c|c|}
		\hline
		 Dimension $d$ & Averaged value & Standard deviation & Relative error (\%) & Theoretical solution\\
		\hline
 5 & 0.37901
 &  0.00312  &  0.97 &  0.382727\\
		\hline
15 & 0.25276
& 0.00235 & 1.17 & 0.255754\\
		\hline
	\end{tabular}
	\caption{Estimate of $u(0, x_0=1_{5})$ on the Monge Ampere problem \eqref{Monge} with  $N = 120$.  
	Average and standard deviation observed over 10 independent runs are reported.}
	\label{fig: table results Monge Ampere}
\end{table}

\subsection{Portfolio selection}


We consider a portfolio selection problem formulated as follows. There are $n$ risky assets of uncorre\-lated price process $P$ $=$ $(P^1,\ldots,P^n)$ with dynamics 
\begin{align}
\di P_t^i  & = \; P_t^i  \sigma(V_t^i)   \big[ \lambda_i (V_t^i)  \di t +   \di W_t^i\big], \quad i=1,\ldots,n,
\end{align} 
where $W$ $=$ $(W^1,\ldots,W^n)$  is a $n$-dimensional Brownian motion, $b$ $=$ $(b^1,\ldots,b^n)$ is the rate of return of the assets, 
$\lambda$ $=$ $(\lambda^1,\ldots,\lambda^n)$  is the risk premium of the assets, 
$\sigma$ is a positive function (e.g. $\sigma(v)$ $=$ $e^v$ corresponding to the Scott model),  
and  $V$ $=$ $(V^1,\ldots,V^n)$  is the  volatility factor modeled by an Ornstein-Uhlenbeck (O.U.) process
\begin{align} \label{dynY} 
 \di V_t^i  &= \;  \kappa_i[ \theta_i - V_t^i] \di t + \nu_i   \di B^i_t, \quad i=1,\ldots,n, 
\end{align} 
with  $\kappa_i, \theta_i, \nu_i$ $>$ $0$, and 
$B$ $=$ $(B^1,\ldots,B^n)$ a $n$-dimensional Brownian motion, s.t. $d<W^i,B^j>$ $=$ $\delta_{ij} \rho_{ij} dt$, with $\rho_{i}$ $:=$ $\rho_{ii}$ $\in$ $(-1,1)$.  An agent can invest at any time an amount $\alpha_t$ $=$ $(\alpha_t^1,\ldots,\alpha_t^n)$  in the stocks, 
which generates a wealth process $\Xc$ $=$ $\Xc^\alpha$ governed by 
\begin{align}
\di \Xc_t  &= \;   \sum_{i=1}^n  \alpha_t^i  \sigma(V_t^i) \big[ \lambda_i(V_t^i) \di t +   \di W_t^i \big].  
\end{align}
The objective of the agent is to  maximize her expected utility from terminal wealth: 
\begin{align}
\E \big[ U(\Xc_T^\alpha) ] \quad  & \leftarrow  \quad \mbox{maximize over } \alpha 
\end{align}
It is well-known that the solution to this problem can be characterized by the dynamic programming method (see e.g. \cite{pha09}), which leads to the Hamilton-Jacobi-Bellman for the value function on 
$[0,T)\times\R\times\R^n$:  
\begin{equation*}
    \begin{cases}
                \partial_t u  + \Sum_{i=1}^n \big[ \kappa_i(\theta_i - v_i) \partial_{v_i} u + \frac{1}{2} \nu_i^2  \partial_{v_i}^2  u  \big] \; = \; 
       \frac{1}{2} R(v)  \frac{(\partial_\mrx u)^2}{  \partial_{\mrx\mrx}^2 u}  + \sum_{i=1}^n \big[  \rho_i \lambda_i(v_i) \nu_i \frac{\partial_\mrx u \partial_{\mrx v_i}^2 u}{ \partial_{\mrx\mrx}^2 u}  + 
        \frac{1}{2}\rho_i^2 \nu_i^2  \frac{(\partial^2_{\mrx v_i} u)^2}{ \partial^2_{\mrx\mrx} u} \big] 
         &
         \\
        u(T,\mrx,v) \; = \;  U(\mrx),   \quad \quad \mrx \in \R,  \;  v \in \R^n, &
    \end{cases}
\end{equation*} 
with a Sharpe ratio $R(v)$ $:=$ $|\lambda(v)|^2$, for $v$ $=$ $(v_1,\ldots,v_n)$ $\in$ $(0,\infty)^n$.  
The optimal portfolio strategy is then  given in feedback form by 
$\alpha_t^*$ $=$ $\hat a(t,\Xc_t^*,V_t)$, where $\hat a$ $=$ $(\hat a_1,\ldots,\hat a_n)$ is given by 
\begin{align}
\hat a_i(t,\mrx,v)  &= \;   - \frac{1}{\sigma(v_i)} \Big(  \lambda_i(v_i)   \frac{\partial_\mrx u}{ \partial_{\mrx\mrx}^2 u} +   \rho_i \nu_i \frac{\partial_{\mrx v_i}^2 u}{ \partial_{\mrx\mrx}^2 u } \Big),  
\quad (t,\mrx,v=(v_1,\ldots,v_n)) \in [0,T)\times\R\times\R^n, 
\end{align}
for $i$ $=$ $1,\ldots,n$.  This Bellman equation is in the form \eqref{eq:PDEInit}  with
\begin{align} \label{Bellmanf} 
f(t,x,y,z,\gamma) &= \;  \Sum_{i=1}^n \big[ \kappa_i(\theta_i - v_i) z_i  + \frac{1}{2} \nu_i^2 \gamma_{ii} \big] - \frac{1}{2} R(v)  \frac{z_0^2}{\gamma_{00}} 
- \sum_{i=1}^n \big[  \rho_i \lambda_i(v_i) \nu_i \frac{z_0  \gamma_{0i}}{\gamma_{00}}  
+  \frac{1}{2}\rho_i^2 \nu_i^2  \frac{(\gamma_{0i})^2}{ \gamma_{00}} \big], 
\end{align}
for  $x$ $=$ $(\mrx,v)$ $\in$ $\R^{n+1}$, $z$ $=$ $(z_0,\ldots,z_n)$ $\in$ $\R^{n+1}$, $\gamma$ $=$ $(\gamma_{ij})_{0\leq i,j\leq  n}$ $\in$ $\S^{n+1}$, 
and  displays a high non-linearity in the Hessian argument $\gamma$.

\vspace{2mm}

The truncation operator indexed by a parameter $p$  is chosen equal to
\begin{align*}
 \Tc_p(X_t^{0,x}) = \min\big\{\max[x + \mu t  -  \sigma \sqrt{t} \phi_p ,  X_t^{0,x} ],  x + \mu t+  \sigma \sqrt{t} \phi_p \big\},
\end{align*}
where $\phi_p$ $=$  $\mathcal{N}^{-1}(p)$,  $\mathcal{N}$ is the CDF of a unit centered Gaussian random variable. 
We use neural networks with 2 hidden layers and  $d+10$ neurons per layer. 
We shall test this example when the utility function $U$ is of exponential form: $U(x)$ $=$ $-\exp(-\eta x)$, with $\eta$ $>$ $0$,  
and under different cases for which  closed-form solutions are available: 
\begin{itemize}
\item[(1)]  {\it Merton problem.} This corresponds to a degenerate case where the factor $V$, hence the volatility $\sigma$ and  the risk premium $\lambda$   are  constant, 
so that the generator of  Bellman equation reduces to 
\begin{equation}\label{PDE MERTON}
f(t,x,y,z,\gamma) \; = \;  - \frac{1}{2} |\lambda|^2  \frac{z^2}{\gamma}, \quad \;\;  (t,x,y,z) \in [0,T]\times\R\times\R\times\R, 
\end{equation} 
with explicit solution given by: 
\begin{align}
u(t,x) & = \;  e^{-(T-t) \frac{|\lambda|^2}{2}}  U(x), \quad \hat a_i \; = \; \frac{\lambda_i}{\eta\sigma}. 
\end{align}
We train with the forward process
\begin{align}
X_{k+1} & = \; X_k + \lambda  \Delta t_k + \Delta W_k, \quad k=0,\ldots,N, \; X_0 \;  = \; x_0.
\end{align}

\item[(2)]  {\it One risky asset: $n$ $=$ $1$}. A quasi-explicit solution is provided in \cite{zari01}: 
\begin{align}
u(t,\mrx,v) & = \; U(\mrx) w(t,v), \quad \mbox{ with }  w(t,v) \; = \; \Big\| \exp\Big(- \frac{1}{2} \int_t^T R(\hat V_s^{t,v}) ds \Big) \Big\|_{_{L^{1-\rho^2}}} 
\end{align}
where $\hat V_s^{t,v}$ is the solution  to the modified O.U.  model
\begin{align}
 \di\hat V_s  &= \;  \big[ \kappa (\theta  -  \hat V_s) -   \rho \nu  \lambda(\hat V_s)  \big] \di s + \nu   \di B_s, \quad s \geq t, \; \hat V_t = v.   
\end{align}
We test our algorithm with $\lambda(v)$ $=$ $\lambda v$, $\lambda$ $>$ $0$, for which we have an explicit solution: 
\begin{align}
w(t,v) & = \; \exp\big( - \phi(t) \frac{v^2}{2}  - \psi(t)v -  \chi(t) \big),   \quad (t,v) \in [0,T]\times\R, 
\end{align}
where $(\phi,\psi,\chi)$ are solutions of the Riccati system of ODEs: 
\begin{align}
\dot \phi - 2 \bar\kappa \phi  - \nu^2(1-\rho^2) \phi^2  + \lambda^2 & = \;0, \quad \phi(T) \;  = \; 0,  \\
\dot\psi - (\bar\kappa + \nu^2(1-\rho^2)\phi)\psi  +  \kappa \theta \phi & = \; 0, \quad \psi(T) \; = \; 0, \\
\dot\chi +  \kappa\theta\psi -  \frac{\nu^2}{2} ( - \phi + (1-\rho^2)\psi^2) & = \; 0, \quad \chi(T) \; = \; 0,  
\end{align}
with $\bar\kappa$ $=$ $\kappa + \rho \nu  \lambda$,  and explicitly given by (see e.g. Appendix in \cite{schzhu99}) 
\begin{align}
\phi(t) &= \; \lambda^2 \frac{\sinh(\hat\kappa(T-t))}{\hat\kappa\cosh(\hat\kappa(T-t)) + \bar\kappa \sinh(\hat\kappa(T-t))} \\
\psi(t) &= \;  \lambda^2 \frac{\kappa\theta}{\hat\kappa} \frac{\cosh(\hat\kappa(T-t)) - 1} {\hat\kappa\cosh(\hat\kappa(T-t)) + \bar\kappa \sinh(\hat\kappa(T-t))} \\
\chi(t) & = \; \frac{1}{2(1-\rho^2)} \ln \big[ \cosh(\hat\kappa(T-t)) + \frac{\bar\kappa}{\hat\kappa} \sinh(\hat\kappa(T-t)) \big]  - \frac{1}{2(1-\rho^2)} \bar\kappa(T-t) \\
& \;\;\; - \; \lambda^2 \frac{(\kappa\theta)^2}{\hat\kappa^2} \Big[  \frac{\sinh(\hat\kappa(T-t))} {\hat\kappa\cosh(\hat\kappa(T-t)) + \bar\kappa \sinh(\hat\kappa(T-t))}  - (T-t) \Big] \\
& \;\;\; - \;  \lambda^2 \frac{(\kappa\theta)^2 \bar\kappa}{\hat\kappa^3}   \frac{\cosh(\hat\kappa(T-t)) - 1} {\hat\kappa\cosh(\hat\kappa(T-t)) + \bar\kappa \sinh(\hat\kappa(T-t))},   
\end{align}
with $\hat\kappa$ $=$ $\sqrt{\kappa^2 + 2 \rho\nu\lambda\kappa +  \gamma^2\lambda^2}$. 
We  train with the forward process
\begin{align}
\Xc_{k+1} & = \; \Xc_k +  \lambda \theta \Delta t_k + \Delta W_k,   \quad k=0,\ldots,N-1, \;\;  \Xc_0 =  \mrx_0, \\
V_{k+1} & = \; V_k  + \nu \Delta B_k, \;\;\;   \quad k=0,\ldots,N-1, \;\; V_0 \;  = \; \theta.
\end{align}

\item[(3)] {\it No leverage effect, i.e.,  $\rho_i$ $=$ $0$, $i$ $=$ $1,\ldots,n$}.  In this case, there is a quasi-explicit solution given by 
\begin{align}\label{PDE : No Leverage}
u(t,\mrx,v) & = \; U(\mrx) w(t,v), \;\;  \mbox{ with }  w(t,v) \; = \; \E \Big[ \exp\Big(- \frac{1}{2} \int_t^T R(V_s^{t,v}) ds \Big) \Big], \;\; (t,v) \in [0,T]\times\R^n, 
\end{align}
where $V^{t,v}$ is the solution to \eqref{dynY}, starting from $v$ at time $t$. 
We test our algorithm with $\lambda_i(v)$ $=$ $\lambda_i v_i$, $\lambda_i$ $>$ $0$, $i$ $=$ $1,\ldots,n$, $v$ $=$ $(v_1,\ldots,v_n)$, 
for which we have an explicit solution given by 
\begin{align}
w(t,v) & = \; \exp\Big( - \sum_{i=1}^n \big[  \phi_i(t) \frac{v_i^2}{2}  + \psi_i(t)v_i + \chi_i(t) \big] \Big),   \quad (t,v) \in [0,T]\times\R^n, \\
\phi_i(t) &= \; \lambda_i^2 \frac{\sinh(\hat\kappa_i(T-t))}{\kappa_i \sinh(\hat\kappa_i(T-t)) + \hat\kappa_i\cosh(\hat\kappa_i(T-t))} \\
\psi_i(t) &= \;  \lambda_i^2 \frac{\kappa_i\theta_i}{\hat\kappa_i} \frac{\cosh(\hat\kappa_i(T-t)) - 1} {\kappa_i \sinh(\hat\kappa_i(T-t)) + \hat\kappa_i\cosh(\hat\kappa_i(T-t))} \\
\chi_i(t) & = \; \frac{1}{2} \ln \big[ \cosh(\hat\kappa_i(T-t)) + \frac{\kappa_i}{\hat\kappa_i} \sinh(\hat\kappa_i(T-t)) \big]  - \frac{1}{2} \kappa_i(T-t) \\
& \;\;\; - \; \lambda_i^2 \frac{(\kappa_i\theta_i)^2}{\hat\kappa_i^2} \Big[  \frac{\sinh(\hat\kappa_i(T-t))} {\hat\kappa_i\cosh(\hat\kappa_i(T-t)) + \kappa_i \sinh(\hat\kappa_i(T-t))}  - (T-t) \Big] \\
& \;\;\; - \;  \lambda^2 \frac{(\kappa_i\theta_i)^2 \kappa_i}{\hat\kappa_i^3}   \frac{\cosh(\hat\kappa_i(T-t)) - 1} {\hat\kappa_i\cosh(\hat\kappa_i(T-t)) + \kappa_i \sinh(\hat\kappa_i(T-t))},   
\end{align}
with $\hat\kappa_i$ $=$ $\sqrt{\kappa_i^2 +  \nu_i^2\lambda_i^2}$. 
We train with the forward process 
\begin{align}
\Xc_{k+1} & = \; \Xc_k  + \Delta W_k,   \quad k=0,\ldots,N-1, \;\;  \Xc_0 \; = \;  \mrx_0, \\
V_{k+1}^i & = \; V_k^i  +  \nu_i \Delta B_k^i, \;\;\;   \quad k=0,\ldots,N-1, \;\; V_0^i \;  = \; \theta_i,
\end{align}
with $<W,B^i>_t$ $=$ $0$. 
\end{itemize}

\vspace{3mm}

\textbf{Merton Problem.} 
We take $\eta = 0.5$, $\lambda = 0.6$, $T$ $=$ $1$, $N = 120$,  and $\sigma(v) = e^v$. 
We plot the neural networks approximation of $u, D_x u, D^2_x u, \alpha$ (in blue) together with their analytic values (in orange). For comparison with Figures \ref{fig: table results Merton Jentzen} and \ref{fig: table results Merton ours},  we report the error on the gradient and the initial control.
In practice, after empirical tests,  we choose $p= 0.98$ for the truncation. 

\begin{table}[H]
	\centering
	\begin{tabular}{|c|c|c|c|c|c|}
		\hline
		 & Averaged value & Standard deviation & Theoretical value & Relative error (\%)\\
		\hline
	$u(0, x_0=1)$	&  -0.50561
 & 0.00029 & -0.50662 & 0.20\\
		\hline
	$D_x	u(0, x_0=1)$ & 0.25081
 & 0.00088 &  0.25331 & 0.99\\
		\hline
	$\alpha(0, x_0=1)$ & 0.83552
 & 0.02371 & 0.80438 & 3.87 \\
		\hline
	\end{tabular}
	\caption{Estimate of the solution, its derivative and the optimal control at the initial time $t=0$ in the Merton problem \eqref{PDE MERTON}. 
	Average and standard deviation observed over 10 independent runs are reported. } 
	\label{fig: table results Merton}
\end{table}

\begin{figure}[h!]
\begin{minipage}[b]{0.32\linewidth}
  \centering
 \includegraphics[width=\textwidth]{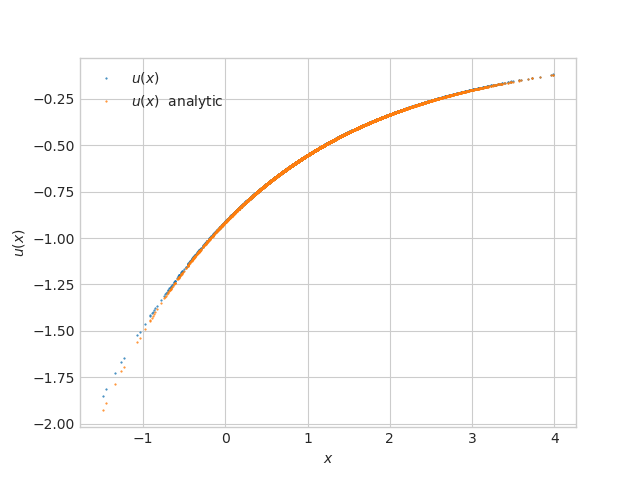}
 \caption*{$Y$ at date $t=0.5042$.}
 \end{minipage}
\begin{minipage}[b]{0.32\linewidth}
  \centering
 \includegraphics[width=\textwidth]{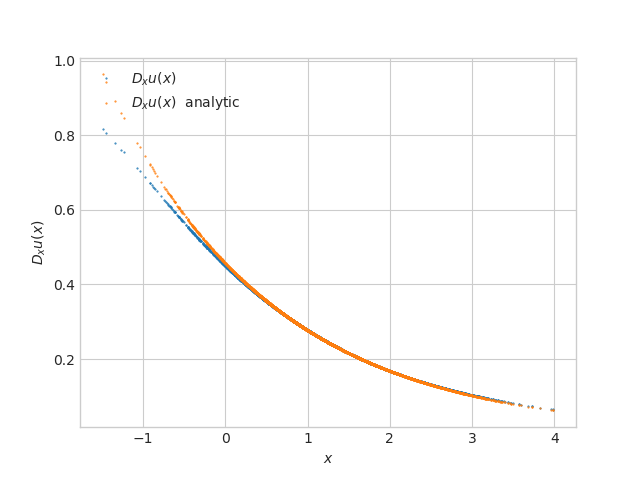}
 \caption*{$Z$ at date $t=0.5042$ }
 \end{minipage}
 \begin{minipage}[b]{0.32\linewidth}
  \centering
 \includegraphics[width=\textwidth]{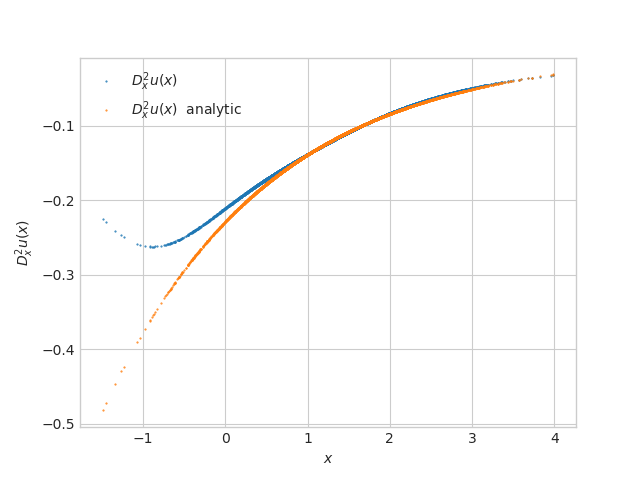}
 \caption*{$\Gamma$ at date $t=0.5042$ }
 \end{minipage}
 \begin{minipage}[b]{0.32\linewidth}
  \centering
 \includegraphics[width=\textwidth]{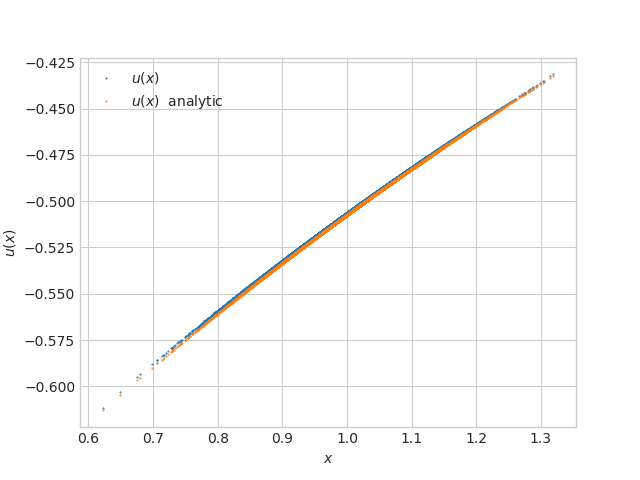}
 \caption*{$Y$ at date $t=0.0084$.}
 \end{minipage}
\begin{minipage}[b]{0.32\linewidth}
  \centering
 \includegraphics[width=\textwidth]{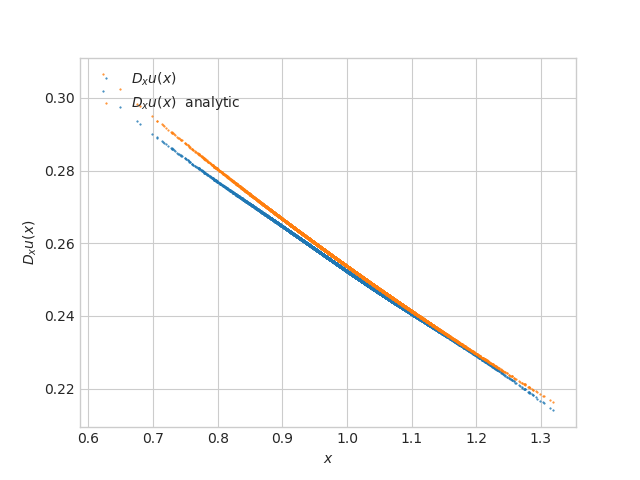}
 \caption*{$Z$ at date $t=0.0084$ }
 \end{minipage}
 \begin{minipage}[b]{0.32\linewidth}
  \centering
 \includegraphics[width=\textwidth]{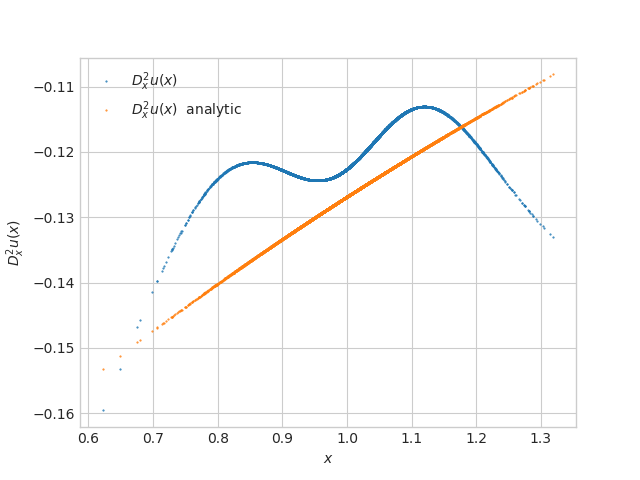}
 \caption*{$\Gamma$ at date $t=0.0084$ }
 \end{minipage}
 \caption{\label{fig: Merton}  Estimates of the solution and its derivatives on the Merton problem \eqref{PDE MERTON}  using $120$ time steps.}
 \end{figure}
 
\begin{figure}[H]
   \begin{minipage}[c]{.49\linewidth}
          \includegraphics[width=0.9\linewidth]{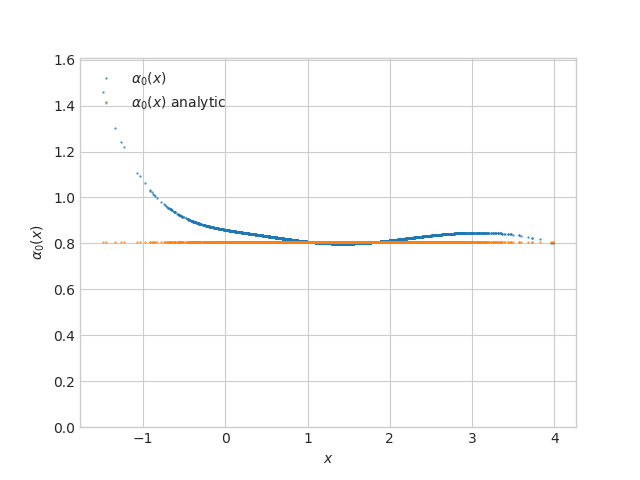}
           \caption*{$\alpha$ at date $t=0.5042$.}
   \end{minipage} \hfill
   \begin{minipage}[c]{.49\linewidth}
      \includegraphics[width=0.9\linewidth]{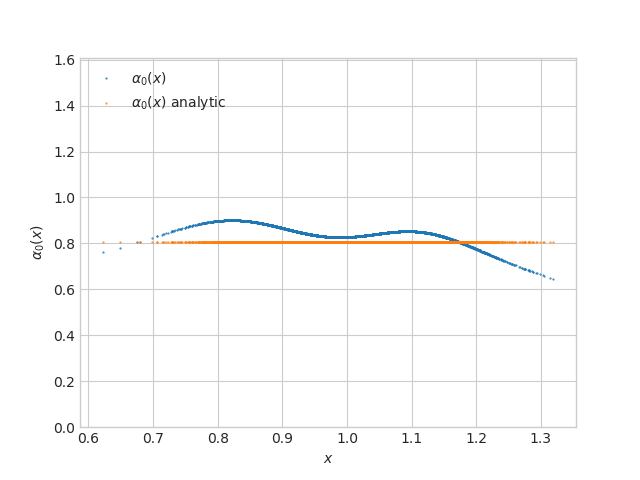}
       \caption*{$\alpha$ at date $t=0.0084$.}
   \end{minipage}
   \caption{Estimates of the optimal control $\alpha$ on the Merton problem \eqref{PDE MERTON}.}
   \label{fig: control Merton}
\end{figure}

\textbf{One asset ($n$ $=$ $1$) in Scott volatility model.}
We take $\eta = 0.5$, $\lambda = 1.5$, $\theta = 0.4$, $\nu = 0.4$, $\kappa = 1$, $\rho = -0.7 $. For all tests we choose $T$ $=$ $1$, $N = 120$,  and $\sigma(v) = e^v$. In practice, after empirical tests,  we choose $p= 0.98$ for the truncation.

\begin{table}[H]
	\centering
	\begin{tabular}{|c|c|c|c|}
		\hline
		 Averaged value & Standard deviation & Relative error (\%)\\
		\hline
		-0.53431 & 0.00070  & 0.34\\
		\hline
	\end{tabular}       
	\caption{Estimate of $u(0, \mrx_0=1, \theta)$ on the One Asset problem with stochastic volatility  ($d=2$). Average and standard deviation observed over 10 independent runs are reported. The exact solution is $-0.53609477$.}
	\label{fig: table results One asset d = 2}
\end{table}

\vspace{1mm}

\textbf{No Leverage in Scott model.} 
In the case with one asset ($n$ $=$ $1$), we take $\eta = 0.5$, $\lambda = 1.5$, $\theta = 0.4$, $\nu = 0.2$, $\kappa = 1$. For all tests we choose $T$ $=$ $1$, $N = 120$, and $\sigma(v) = e^v$. In practice, after empirical tests,  we choose $p= 0.95$ for the truncation.

\begin{table}[H]
    \centering
    \begin{tabular}{|c|c|c|c|c|}
    \hline
        Dimension $d$ & Averaged value & Standard deviation & Relative error (\%) & Theoretical solution\\
		\hline
2 & -0.49980 & 0.00073 & 0.35 & -0.501566\\
	\hline
5 & -0.43768 & 0.00137 & 0.92 & -0.441765\\
	\hline
8 & -0.38720  & 0.00363 & 1.96 & -0.394938\\
	\hline  
10 &  -0.27920 & 0.05734 & 1.49 & -0.275092\\
	\hline
    \end{tabular}
    \caption{Estimate of $u(0, \mrx_0=1, \theta)$  on the No Leverage problem \eqref{PDE : No Leverage}.  
	Average and standard deviation observed over 10 independent runs are reported.} 
    \label{tab:my_label}
\end{table}


In the case with four assets ($n$ $=$ $4$, $d=5$), we take $\eta = 0.5$, $\lambda = \begin{pmatrix} 1.5 & 1.1 & 2. & 0.8 \end{pmatrix}$, $\theta = \begin{pmatrix} 0.1 & 0.2 & 0.3 & 0.4 \end{pmatrix}$, 
$\nu = \begin{pmatrix}0.2 & 0.15 & 0.25 & 0.31 \end{pmatrix}$, $\kappa = \begin{pmatrix}1. & 0.8 & 1.1 & 1.3\end{pmatrix}$.

	

In the case with seven assets ($n$ $=$ $7$, $d=8$) we take $\eta = 0.5$, $\lambda = \begin{pmatrix} 1.5 & 1.1 & 2. & 0.8 & 0.5 & 1.7 & 0.9 \end{pmatrix}$,\\ $\theta = \begin{pmatrix} 0.1 & 0.2 & 0.3 & 0.4 & 0.25 & 0.15 & 0.18  \end{pmatrix}$, $\nu = \begin{pmatrix}0.2 & 0.15 & 0.25 & 0.31 & 0.4 & 0.35 & 0.22 \end{pmatrix}$,\\ $\kappa = \begin{pmatrix}1. & 0.8 & 1.1 & 1.3 & 0.95 & 0.99 & 1.02 \end{pmatrix}$.


In the case with nine assets ($n$ $=$ $9$, $d=10$), we take $\eta = 0.5$, $\lambda = \begin{pmatrix} 1.5 & 1.1 & 2. & 0.8 & 0.5 & 1.7 & 0.9 & 1. & 0.9\end{pmatrix}$, $\theta = \begin{pmatrix} 0.1 & 0.2 & 0.3 & 0.4 & 0.25 & 0.15 & 0.18 & 0.08 & 0.91 \end{pmatrix}$, $\nu = \begin{pmatrix}0.2 & 0.15 & 0.25 & 0.31 & 0.4 & 0.35 & 0.22 & 0.4 & 0.15\end{pmatrix}$, $\kappa = \begin{pmatrix}1. & 0.8 & 1.1 & 1.3 & 0.95 & 0.99 & 1.02 & 1.06 & 1.6\end{pmatrix}$.

		

\vspace{2mm}

Hamilton-Jacobi-Bellman equation from portfolio optimization is a typical example of full-nonlinearity in the second order derivative, and the above results show that our algorithm performs quite 
well  up to dimension $d$ $=$ $8$, but gives a high variance in dimension $d$ $=$ $10$. 
\vspace{2mm}

\textbf{Comparison with an implicit version of the scheme.} 
As explained in Remark \ref{rem: implicit explicit}, an alternative option for the estimation of the Hessian is to approximate it by the automatic differentiation of the current neural network for the $Z$ component. It corresponds to the replacement of $ D \hat \Zc_{i+1}(\Tc(X_{t_{i+1}}))$ by $D \Zc_i(\Tc(X_{t_i}));\theta) $ in \eqref{eq:scheme}. An additional change has to be made to the method for it to work. At the last optimization step (for time step $t_0 = 0$), we notice empirically that the variable $\Gamma_0$ is not able to properly learn the initial Hessian value at all. Therefore for this last step we use variables $Y_0,Z_0$ and an explicit estimation of the second order derivative given by $ D \hat \Zc_{1}(\Tc(X_{t_{1}}))$. We see in Table \ref{fig: table results Merton Implicit} that the results for the Merton problem are very similar to the ones from Table \ref{fig: table results Merton} for the splitting scheme but with a worse estimation of the Hessian and optimal control (the error is multiplied by around 1.5).  When we tested this implicit scheme on the Monge Ampere problem we also faced computational problems during the optimization step of Tensorflow. The numerical computation of the gradient of the objective function for the backpropagation step, more precisely for the determinant part, often gives rise to matrix invertibility errors which stops the algorithm execution. For these two reasons, we focused our study on the explicit scheme.

\begin{table}[H]
	\centering
	\begin{tabular}{llllll}
		\hline\noalign{\smallskip}
		 & Average & Std & True value & Relative error (\%)\\
		\noalign{\smallskip}\hline\noalign{\smallskip}
	$u(0, x_0=1)$	& -0.50572  &  0.00034
  & -0.50662 & 0.18\\
		
	$D_x	u(0, x_0=1)$ & 0.25091
 & 0.00067 &  0.25331 & 0.95 \\
		
	$\alpha(0, x_0=1)$ & 0.85254
 & 0.01956 & 0.80438 & 5.99 \\
		\noalign{\smallskip}\hline
	\end{tabular}
	\caption{Estimate of the solution, its derivative and the optimal control at the initial time $t=0$ in the Merton problem \eqref{PDE MERTON} with implicit estimation of the Hessian. 
	Average and standard deviation (Std) observed over 10 independent runs are reported}  
	\label{fig: table results Merton Implicit}
\end{table}

\vspace{2mm}

\begin{figure}[H]
    \begin{minipage}[c]{.49\linewidth}
          \includegraphics[width=0.9\linewidth]{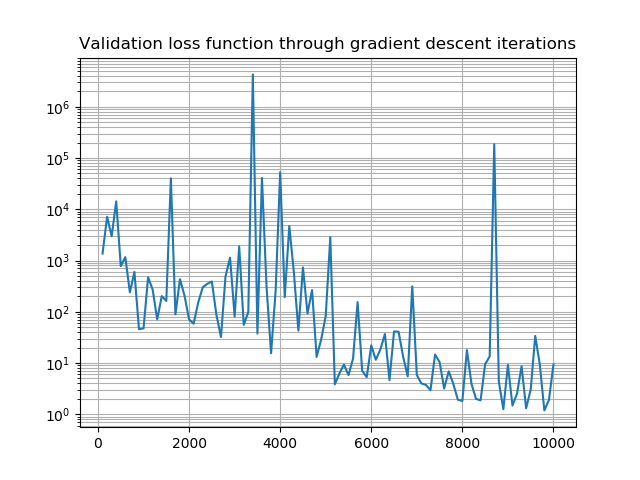}
   \end{minipage} \hfill
   \begin{minipage}[c]{.49\linewidth}
      \includegraphics[width=0.9\linewidth]{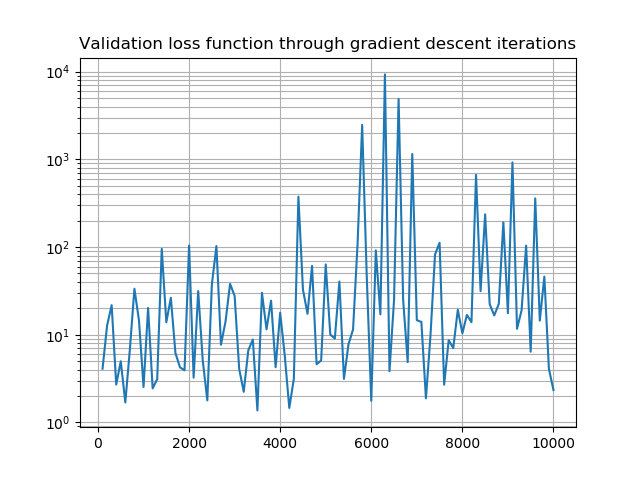}
   \end{minipage}
    \caption{Learning curve in logarithmic scale for the scheme \cite{BEJ19} on the Merton problem \eqref{PDE MERTON} with $N$ $=$ $20$ times steps on the left and $N$ $=$ $120$ 
    time steps on the right. The maturity is $T=1$}. 10000 gradient descent iterations were conducted.
    \label{fig: BEJ}
\end{figure}

\textbf{Comparison with the 2BSDE scheme of \cite{BEJ19}.} 
We  conclude this paper with a comparison of our algorithm with the global scheme of \cite{BEJ19}, called Deep 2BDSE. The tests below concern the Merton problem \eqref{PDE MERTON} but similar behavior happens on the other examples with stochastic volatilities. This scheme 
was implemented in the original paper only for  small number of time steps (e.g. $N$ $=$ $30$). Thus we tested this algorithm on two discretizations, respectively with $N$ $=$ $20$ and $N$ $=$ $120$ time steps, as shown in Figure \ref{fig: BEJ}, for $T=1$ where we plotted the learning curve of the Deep BSDE method.  These curves correspond to the values taken by the loss function  during the gradient descent iterations. For this algorithm the loss function to minimize in the training of neural networks is defined as the mean $L^2$ error between the generated $Y_N$ value and the true terminal condition $g(X_N)$. We observe that for this choice of maturity $T=1$ the loss function oscillates during the training process and does not vanish. As a consequence the Deep 2BSDE does not converge in this case.  Even when decreasing the learning rate, we noticed that we cannot obtain the convergence of the scheme.  
\paragraph{}
However, the Deep 2BSDE method does converge for small maturities $T$, as illustrated in Table \ref{fig: table results Merton Jentzen}  with $T=0.1$ and different values for the number of time steps $N$.  Nevertheless, even if the value function is well approximated, the estimation of the gradient and control did not converge  (the corresponding variance is very large), in comparison with our scheme whereas the gradient is very well approximated and the control is quite precise. We also have a much smaller variance in the results. Table  \ref{fig: table results Merton ours} shows the results obtained by our method with $T=0.1$ in order to compare it with the performance of \cite{BEJ19}. It illustrates the limitations of the global approach and justifies our introduction of a local method.

\begin{table}[H]
	\centering
	\begin{tabular}{|c|c|c|c|c|c|}
		\hline
		 & $N$ & Averaged value & Standard deviation & Theoretical value & Relative error (\%)\\
		\hline
$u(0, x_0= 1)$	&	5 & -0.60667 & \textbf{0.01588 }
 & -0.59571 & 1.84 \\
		\hline
$u(0, x_0= 1)$	&	10 &  -0.59841 &   0.02892
 & -0.59571 & 0.45  \\
		\hline
$u(0, x_0= 1)$	&	20 &  \textbf{-0.59316} & 0.04251
 & -0.59571 & \textbf{0.43} \\
		\hline\hline
		$D_x u(0, x_0= 1)$	&	5 &  \textbf{0.09668} & \textbf{0.25630}
 &  0.29786 & \textbf{67.54} \\
		\hline
$D_x u(0, x_0= 1)$	&	10 &  0.03810 &  0.44570
 &  0.29786 & 93.36  \\
		\hline
$D_x u(0, x_0= 1)$	&	20 & 0.07557 & 0.55030
 &  0.29786 & 74.63  \\
	\hline\hline
	$\alpha(0, x_0=1)$ & 5 & -0.15243
 & 0.61096 & 0.80438 & 118.95 \\
		\hline
	$\alpha(0, x_0=1)$ & 10  & \textbf{0.59971}
 & 1.97906 & 0.80438  & \textbf{25.44}\\
		\hline
	$\alpha(0, x_0=1)$ & 20 & 0.28385
 & \textbf{0.43775} & 0.80438 & 64.71\\
		\hline
	\end{tabular}
	\caption{Estimate of the solution, its derivative and the optimal control at the initial time $t=0$ in the Merton problem \eqref{PDE MERTON} with maturity $T = 0.1$ for the \cite{BEJ19} scheme. 
	Average and standard deviation observed over 10 independent runs are reported.}
	\label{fig: table results Merton Jentzen}
\end{table}

\begin{table}[H]
	\centering
	\begin{tabular}{|c|c|c|c|c|c|}
		\hline
		 & $N$ & Averaged value & Standard deviation & Theoretical value & Relative error (\%)\\
		\hline
$u(0, x_0= 1)$	&	5 &  \textbf{-0.59564} &  0.01136
 & -0.59571 & \textbf{0.01} \\
		\hline
$u(0, x_0= 1)$	&	10 & -0.59550 &  \textbf{0.00037} 
 & -0.59571 & 0.04  \\
		\hline
$u(0, x_0= 1)$	&	20 & -0.59544 & 0.00054
 & -0.59571 & 0.04 \\
		\hline \hline
		$D_x u(0, x_0= 1)$	&	5 &  0.29848 &  \textbf{0.00044}
 &  0.29786 & 0.21 \\
		\hline
$D_x u(0, x_0= 1)$	&	10 &  0.29842
  &  0.00084
 &  0.29786 &  0.19 \\
		\hline
$D_x u(0, x_0= 1)$	&	20 & \textbf{0.29785}  & 0.00054
 &  0.29786 &  \textbf{0.001} \\
	\hline\hline
	$\alpha(0, x_0=1)$ & 5 & \textbf{0.82322}
 & \textbf{0.01014} & 0.80438 & \textbf{2.34}\\
		\hline
	$\alpha(0, x_0=1)$ & 10  & 0.85284
 & 0.07565 & 0.80438 & 6.02 \\
		\hline
	$\alpha(0, x_0=1)$ & 20 & 0.84201
 & 0.09892 & 0.80438 & 4.68\\
		\hline
	\end{tabular}
	\caption{Estimate of the solution, its derivative and the optimal control at the initial time $t=0$ in the Merton problem \eqref{PDE MERTON} with maturity $T = 0.1$ for our scheme. 
	Average and standard deviation observed over 10 independent runs are reported.}
	\label{fig: table results Merton ours}
\end{table}

\printbibliography

\end{document}